\documentclass[a4paper]{amsart}
\usepackage[T1]{fontenc}
\usepackage[latin1]{inputenc}
\usepackage{amssymb}
\usepackage[all]{xy}
\usepackage{enumerate}
\usepackage[lite]{amsrefs}

\newcommand*{\abs}[1]{\lvert#1\rvert}
\newcommand*{\norm}[1]{\lVert#1\rVert}
\newcommand*{\enorm}[1]{\lvert\lVert#1\rVert\rvert}

\newcommand*{\cl}[1]{\overline{#1}}
\newcommand*{\blank}{{\llcorner\!\!\lrcorner}}

\newcommand*{\Born}{\mathsf{Bor}}
\newcommand*{\Mod}{\mathsf{Mod}}
\newcommand*{\Ho}{\mathsf{Ho}}
\newcommand*{\Der}{\mathsf{Der}}

\newcommand*{\ima}{\mathrm{i}}
\newcommand*{\ID}{\mathrm{id}}

\newcommand*{\Mult}{\mathcal{M}}
\newcommand*{\Sch}{\mathcal{S}}
\newcommand*{\CCINF}{\mathcal{D}}
\newcommand*{\CINF}{\mathcal{E}}

\newcommand*{\hot}{\mathbin{\hat{\otimes}}}
\newcommand*{\Lhot}{\mathbin{\hat{\otimes}}^{\mathbb{L}}}
\newcommand*{\cross}{\mathbin{\ltimes}}
\newcommand*{\defeq}{\mathrel{:=}}
\newcommand*{\into}{\rightarrowtail}
\newcommand*{\prto}{\twoheadrightarrow}
\newcommand*{\congto}{\overset{\cong}\to}
\newcommand*{\inOb}{\in_O}   

\DeclareMathOperator{\coker}{coker}
\DeclareMathOperator{\Hom}{Hom}
\DeclareMathOperator{\End}{End}
\DeclareMathOperator{\Aut}{Aut}

\DeclareMathOperator{\Right}{\mathbb{R}}
\DeclareMathOperator{\Left}{{\mathbb{L}}}
\DeclareMathOperator{\Ext}{Ext}
\DeclareMathOperator{\Tor}{Tor}
\DeclareMathOperator{\cone}{cone}
\DeclareMathOperator{\Smooth}{Ess}
\DeclareMathOperator{\Rough}{Rgh}
\DeclareMathOperator{\Gl}{Gl}
\DeclareMathOperator{\ad}{ad}

\newcommand*{\op}{\mathrm{op}}
\newcommand*{\dR}{\mathrm{dR}}

\newcommand*{\T}{{\mathcal{T}}}
\newcommand*{\LG}{{\mathfrak{g}}}
\newcommand*{\LGP}{{\mathfrak{p}}}
\newcommand*{\UE}{{\mathrm{U}}}
\newcommand*{\Cat}{{\mathfrak{C}}}
\newcommand*{\Fil}{{\mathcal{F}}}
\newcommand*{\Exact}{{\mathsf{Exact}}}
\newcommand*{\Subtri}{{\mathfrak{T}}}

\newcommand*{\Pol}{{\mathcal{P}}}
\newcommand*{\Hol}{{\mathcal{O}}}

\newcommand*{\C}{{\mathbb{C}}}
\newcommand*{\F}{{\mathbb{F}}}
\newcommand*{\R}{{\mathbb{R}}}
\newcommand*{\Z}{{\mathbb{Z}}}

\newcommand*{\N}{{\mathbb{N}}}
\newcommand*{\Q}{{\mathbb{Q}}}
\newcommand*{\Torus}{{\mathbb{T}}}

\newcommand*{\brd}{-\hspace{0pt}}
\newcommand*{\nbd}{\nobreakdash-\hspace{0pt}}

\theoremstyle{plain}
\newtheorem{theorem}{Theorem}
\newtheorem{proposition}[theorem]{Proposition}
\newtheorem{lemma}[theorem]{Lemma}
\newtheorem{corollary}[theorem]{Corollary}

\theoremstyle{definition}
\newtheorem{definition}[theorem]{Definition}
\theoremstyle{remark}

\newtheorem{example}[theorem]{Example}

\begin{document}

\title{Embeddings of derived categories of bornological modules}
\author{Ralf Meyer}
\address{Mathematisches Institut\\
         Westfälische Wilhelms-Universität Münster\\
         Einsteinstr.\ 62\\
         48149 Münster\\
         Germany
}
\email{rameyer@math.uni-muenster.de}

\subjclass[2000]{Primary 46M18; Secondary 18A40, 18E30, 19D55, 46A17, 46M15}

\thanks{This research was supported by the EU-Network \emph{Quantum
    Spaces and Noncommutative Geometry} (Contract HPRN-CT-2002-00280)
  and the \emph{Deutsche Forschungsgemeinschaft} (SFB 478).}

\begin{abstract}
  Let~$A$ be an algebra with a countable basis and let~$B$ be, say, a Fréchet
  algebra that contains~$A$ as a dense subalgebra.  The embedding $A\to B$
  induces a functor from the derived category of $B$\nobreakdash-modules to
  the derived category of $A$\nobreakdash-modules.  In many important
  examples, it happens that this functor is fully faithful.  We study this
  property in some detail, giving several equivalent conditions, examples, and
  applications.
  
  To prepare for this, we explain carefully how to do homological algebra with
  modules over bornological algebras.  We construct the derived category of
  bornological left $A$\nobreakdash-modules and some standard derived
  functors, with special emphasis on the adjoint associativity between the
  tensor product and the internal $\Hom$ functor.  We also discuss the
  category of essential modules over a non-unital algebra and its
  functoriality.
\end{abstract}
\maketitle

\section{Introduction}
\label{sec:intro}

We first consider a motivating example.  Let~$\Torus^2_\theta$ be the
noncommutative $2$\nbd{}torus with parameter $\theta\in\R$.  Functions
on~$\Torus^2_\theta$ have Fourier expansions of the form $\sum_{m,n\in\Z}
a_{mn} U^mV^n$, where $U$ and~$V$ are invertible and satisfy the relation
$$
UV=\exp(2\pi\ima\theta)VU.
$$
The algebra $\Pol(\Torus^2_\theta)$ of \emph{polynomial functions} on
$\Torus^2_\theta$ is defined by the requirement $a_{mn}\in\C[\Z^2]$.  The
algebra $\Sch(\Torus^2_\theta)$ of \emph{smooth functions} on
$\Torus^2_\theta$ is defined by the condition $(a_{mn})\in\Sch(\Z^2)$, that
is, $a_{mn}=O(m+n+1)^{-k}$ for all $k\in\N$.  If $\theta=0$, then
$\Pol(\Torus^2_\theta)$ is the algebra of Laurent series, that is, polynomial
functions on $(\C^\times)^2$, and $\Sch(\Torus^2_\theta)$ is isomorphic to the
algebra of smooth functions on~$\Torus^2$.

Alain Connes has computed the Hochschild and cyclic homology of
$\Sch(\Torus^2_\theta)$ in~\cite{Connes:Noncommutative_Diffgeo}.  He uses a
small free $\Sch(\Torus^2_\theta)$\brd{}bimodule resolution of
$\Sch(\Torus^2_\theta)$.  Such a resolution is rather easy to find for the
dense subalgebra $\Pol(\Torus^2_\theta)$.  The crucial point in Connes's
computation is that the same type of complex still works for
$\Sch(\Torus^2_\theta)$.  There are situations that are similar at first
sight, but where such a resolution does not exist.  For instance, there is no
small free $\ell_1(\Z)$\brd{}bimodule resolution for the convolution algebra
$\ell_1(\Z)$, although there is an evident one for the dense subalgebra
$\C[\Z]$.  This is why the Hochschild homology of $\ell_1(\Z)$ is so hard to
compute.

The above example generalises as follows.  Let $f\colon A\to B$ be a bounded
unital homomorphism between two bornological unital algebras.  If you are
unfamiliar with bornologies, you may think of topological algebras here for
the time being.  Choose any resolution of~$A$ by free bornological
$A$\nbd{}bimodules,
$$
\dotsb \to P_3 \to P_2 \to P_1 \to P_0 \to A.
$$
``Resolution'' means that the chain complex above has a bounded contracting
homotopy.  Now we consider the chain complex $B\hot_A P_\bullet \hot_A B$,
which can be equipped with a canonical augmentation map $B\hot_A P_0 \hot_A
B\to B$.  Here~$\hot_A$ denotes the $A$\nbd{}balanced complete bornological
tensor product, which is defined by the universal property that bounded linear
maps $X\hot_A Y\to Z$ correspond to bounded bilinear maps $b\colon X\times
Y\to Z$ that satisfy $b(xa,y)=b(x,ay)$.

The complex $B\hot_A P_\bullet \hot_A B$ is a good candidate for a free
$B$\nbd{}bimodule resolution of~$B$.  It is always a chain complex of free
bornological $B$\nbd{}bimodules.  If it is a resolution of~$B$ as well, we
call $f\colon A\to B$ \emph{isocohomological}.  This is what happens for the
embedding of $\Pol(\Torus^2_\theta)$ into $\Sch(\Torus^2_\theta)$ and makes
the cyclic homology of $\Sch(\Torus^2_\theta)$ computable.  I introduced this
concept in the context of group convolution algebras
in~\cite{Meyer:Primes_Rep} to facilitate some rather technical computations of
coinvariant spaces.  Of course, the information that~$f$ is isocohomological
is particularly useful if~$A$ has a rather small free bimodule resolution.
Proving that a homomorphism is isocohomological is usually difficult and
requires some geometric insight into the algebras at hand.  This is
particularly apparent in the following situation.

Let~$G$ be a finitely generated discrete group and let~$\ell$ be a word
length function on~$G$.  Let $\C[G]$ be the associated group algebra and let
\begin{equation}  \label{eq:def_Schwartz_space}
  \Sch(G) \defeq
  \Bigl\{ f\colon G\to\C \Bigm|
    \sum_{g\in G} \abs{f(g)} (\ell(g)+1)^k<\infty \quad \forall k\in\N
  \Bigr\}.
\end{equation}
Then $\C[G]$ is an algebra with a countable basis and $\Sch(G)$ is a Fréchet
algebra, containing $\C[G]$ as a dense subalgebra.  The embedding $i\colon
\C[G]\to\Sch(G)$ is isocohomological if and only if the slightly simpler chain
complex $\Sch(G)\hot_{\C[G]} P_\bullet\hot_{\C[G]}\C$ is a resolution of~$\C$,
where~$\C$ is equipped with the trivial representation of~$G$ and $P_\bullet$
is a free $\C[G]$\brd{}bimodule resolution of $\C[G]$.  I show
in~\cite{Meyer:Combable_poly} that the homotopy class of this chain complex
only depends on the large scale geometry of~$G$ and construct a contracting
homotopy in the case where~$G$ is a combable group in the sense
of~\cite{Epstein:Automatic}.  In this article, we will see that it is also
contractible if~$G$ has polynomial growth.  The proof reduces this assertion
to an analogous one for nilpotent Lie groups, which can be checked using chain
complexes of differential forms.

There are finitely generated groups for which the embedding $\C[G]\to\Sch(G)$
is not isocohomological because of the following obstruction: a necessary
condition for $\C[G]\to\Sch(G)$ to be isocohomological is that the rational
cohomology groups $H^n(G;\Q)$ must be finite dimensional for all $n\in\N$.

The main goal of this article is to explore other properties of
isocohomological homomorphisms.  Put in a nutshell, if $f\colon A\to B$ is
isocohomological, then all homological computations with $B$\nbd{}modules can
be reduced to homological computations with $A$\nbd{}modules.  There are three
basic constructions for which we verify this statement: the left derived
functor~$\Lhot_B$ of the balanced tensor product $(M,N)\mapsto M\hot_B N$, the
right derived functor $\Right\Hom_B$ of the functor $(M,N)\mapsto
\Hom_B(M,N)$, and the derived category $\Der(B)$ over the category of
bornological $B$\nbd{}modules.

It is important for us to use total derived functors instead of the satellite
functors like $\Ext^n_B(M,N)$ and $\Tor_n^B(M,N)$ that are often called
derived functors.  The latter can be obtained from the total derived functors
$\Right\Hom_B(M,N)$ and $M\Lhot_B N$ by passing to homology.  Thus statements
about total derived functors imply statements about their satellite functors.
However, the passage to homology forgets the bornology, which is an important
part of the structure.  For instance, the chain complex $B\hot_A P_\bullet
\hot_A B$ is a realisation of $B\Lhot_A B$.  The definition of an
isocohomological homomorphism therefore amounts to the condition that the
natural chain map $B\Lhot_A B\to B\Lhot_B B\cong B$ should be a homotopy
equivalence.  This is more than a statement about $\Tor_n^A(B,B)$.

We will see in Theorem~\ref{the:isocoh} that various conditions involving
$\Lhot_B$, $\Right\Hom_B$, and $\Der(B)$ are equivalent to~$f$ being
isocohomological.  For instance, $f$ is isocohomological if and only if the
functor $f^*\colon \Der(B)\to\Der(A)$ induced by~$f$ is fully faithful.  It is
possible to prove these equivalences by hand, playing around cleverly with
free resolutions.  However, the proof gets more transparent and much shorter
if we use some of the more advanced tools of homological algebra.  There are
also other situations where this machinery is useful.  Therefore, a great part
of this article deals with the application of homological algebra to
categories of bornological modules.

Before I explain how this works, I should discuss why I use bornologies
instead of topologies.  A bornology on a vector space is a collection of
bounded subsets that satisfies certain axioms.  Topological vector spaces
always carry canonical bornologies.  As long as we are only dealing with
Fréchet spaces, bornological and topological analysis tend to be equivalent
(see~\cite{Meyer:Born_Top}).  Already for the simplest possible non-Fréchet
spaces like $\C[G]$ or $\Pol(\Torus^2_\theta)$, working topologically creates
several technical problems.  These problems are artificial because they
disappear when we treat these spaces bornologically.  One of the reasons why
homological algebra with bornological vector spaces is particularly nice is
adjoint associativity: the complete bornological tensor product is left
adjoint to the internal $\Hom$\brd{}functor.  In particular, there is a
canonical bornology on a space of bounded linear maps.  In contrast, on spaces
of continuous linear maps, there are dozens of topologies, which are useful
for certain purposes, but do not have good algebraic properties.  The complete
topological tensor product cannot have any right adjoint because it does not
commute with direct sums.

Adjoint associativity allows us to get assertions about homology from
assertions about cohomology (with coefficients).  In contrast, in categories
of topological vector spaces homology seems to carry more information than
cohomology.  More precisely, let $f\colon A\to B$ be a continuous unital
homomorphism between two topological unital algebras.  The homological
statement $B\Lhot_A B\cong B$ implies easily that the induced functor
$f^*\colon \Der(B)\to\Der(A)$ is fully faithful.  However, the converse seems
to be false.  At least, I have no idea how to prove it in general.

Let $\Mod(A)$ be the category of bornological left modules over some
bornological unital algebra~$A$.  This category is never Abelian.  Jean-Pierre
Schneiders and Fabienne Prosmans have promoted the notion of a quasi-Abelian
category (see \cites{Schneiders:Quasi-Abelian,
  Prosmans-Schneiders:Born_Indlim, Prosmans:Derived_analysis}).  The category
$\Mod(A)$ is indeed a quasi-Abelian category.  However, there are only very
few situations where the resulting derived functors and derived categories are
useful.  The reason is that we usually do not want any nontrivial homological
algebra to happen for $A=\C$.  Moreover, we want free modules to be
projective.  Therefore, we do \emph{relative} homological algebra and only
allow resolutions with a bounded linear section.  This setup can be formalised
easily: $\Mod(A)$ with the class of extensions with a bounded linear section
is an exact category in the sense of Daniel Quillen.  Abelian categories and
quasi-Abelian categories are, of course, special cases of exact categories.
Most constructions that work in Abelian categories still work for exact
categories with very mild extra hypotheses.  Therefore, we shall work in the
framework of exact categories.

We let $\Ho(A)$ be the homotopy category of chain complexes over the category
$\Mod(A)$.  This is a triangulated category, and the exact complexes form a
thick, triangulated subcategory $\Exact(A)\subseteq\Ho(A)$.  This is true for
almost any exact category by~\cite{Neeman:Derived_Exact}.  The derived
category $\Der(A)$ is defined as the localisation of $\Ho(A)$ at $\Exact(A)$.
The exact category $\Mod(A)$ has enough projectives and injectives, namely,
the free modules $A\hot X$ and the cofree modules $\Hom(A,X)$.  This yields a
rather explicit description of $\Der(A)$ and derived functors.  Another
important feature of $\Ho(A)$ and $\Der(A)$ is that these categories are
generated in some sense by appropriate subcategories of $\Mod(A)$ (see
Propositions \ref{pro:ModA_generates_HoA} and \ref{pro:free_generate_DerA}).
This idea goes back to~\cite{Boekstedt-Neeman}.  It is very useful to reduce
assertions about chain complexes to assertions about modules.

There is a canonical bornology on the space $\Hom_A(M,N)$ of bounded
$A$\nbd{}module homomorphisms between two $A$\nbd{}modules $M$ and~$N$.  The
functors $\Hom_A$ and $\hot_A$ satisfy important adjointness relations.  The
most general statement is that there is a natural isomorphism of bornological
vector spaces
$$
\Hom_{B,C}(M\hot_A N,X) \cong \Hom_{A,C}(N,\Hom_B(M,X))
$$
if $M$ is a $B,A$\brd{}bimodule, $N$ is an $A,C$\brd{}bimodule, and~$X$ is
a $B,C$\brd{}bimodule.

It is easy to extend the construction of $\Hom_A(M,N)$ and $M\hot_A N$ to the
case where $M$ and~$N$ are chain complexes instead of modules.  Applying the
resulting functors to projective and/or injective resolutions, we define the
derived functors $\Right\Hom_A(M,N)$ and $M\Lhot_A N$.  These are triangulated
functors
$$
\Right\Hom_A\colon \Der(A)^\op\times\Der(A)\to\Ho,
\qquad
\Lhot_A\colon \Der(A^\op)\times\Der(A) \to \Ho.
$$
Here ${}^\op$ denotes opposite categories and algebras.  They satisfy the
appropriate analogue of adjoint associativity.

It is frequently necessary in noncommutative geometry to consider certain
non-unital algebras like $\CCINF(G)$, the convolution algebra of smooth,
compactly supported functions on a locally compact group~$G$.  Although
non-unital, $\CCINF(G)$ has approximate identities, and $\CCINF(G)$ is
projective both as a left and right module over itself.  We call algebras with
these two properties quasi-unital.  A module over an algebra is called
essential if the module action $A\hot M\to M$ is a bornological quotient map
or, equivalently, $A\hot_A M\cong M$.  We let $\Mod(A)$ be the category of
essential modules.  The category $\Mod(\CCINF(G))$ is of great interest
because it is isomorphic to the category of smooth representations of~$G$ on
bornological vector spaces by~\cite{Meyer:Smooth}.  The case of $\CCINF(G)$ is
typical for general quasi-unital algebras: most arguments for the group case
work in this generality.

Any quasi-unital algebra has a multiplier algebra $\Mult(A)$, and the category
of essential modules over~$A$ embeds as a full subcategory in the category of
unital $\Mult(A)$\brd{}modules.  This embedding also preserves projectives, so
that we get a corresponding fully faithful embedding of derived categories
$\Der(A)\subseteq \Der(\Mult(A))$.  This is not an equivalence of categories
because $\Mult(A)$ is not an essential module over~$A$.  Although it may seem
tempting to work with unital $\Mult(A)$\brd{}modules instead of essential
$A$\nbd{}modules, there are important reasons for not doing so (see
Section~\ref{sec:mult}).

We also discuss the functoriality of $\Mod(A)$ and $\Der(A)$.  A bounded
homomorphism $f\colon A\to\Mult(B)$ allows us to view~$B$ as an
$A$\nbd{}bimodule.  We call~$f$ \emph{essential} or just a \emph{morphism}
if~$B$ is essential as an $A$\nbd{}bimodule.  If $A$ and~$B$ are unital, this
just means that $f(1_A)=1_B$.  A morphism is called \emph{proper} if it is a
bounded map from~$A$ into~$B$.  We use the notation $f\colon A\dashrightarrow
B$ to denote a morphism from~$A$ to~$B$.  We write $f\colon A\to B$ if~$f$ is
proper.  To any morphism, we associate functors
$$
f^*,f^!\colon \Mod(B)\to\Mod(A),
\qquad
f_*,f_!\colon \Mod(A)\to\Mod(B).
$$
The functor $f_*$ is right adjoint to~$f^*$ and $f_!$ is left adjoint
to~$f^!$.  If~$f$ is proper, then $f^*=f^!$.  This is reminiscent of the
adjointness properties between functors in sheaf theory that carry the same
names.  However, beyond this formal point there is no further relationship
between these two situations.

For example, if~$H$ is a closed subgroup of a locally compact group~$G$, then
we get an associated morphism $i_H^G\colon \CCINF(H)\to\CCINF(G)$.  It is
proper if and only if~$H$ is open in~$G$.  The functors $(i_H^G)^*\colon
\Mod(G)\to\Mod(H)$, $(i_H^G)_*\colon \Mod(H)\to\Mod(G)$, and $(i_H^G)_!\colon
\Mod(H)\to\Mod(G)$ correspond to restriction, induction, and compact induction
of representations, up to certain relative modular functions.  There is no
classical analogue of the functor $(i_H^G)^!$.

After introducing the above structure on module categories, we turn to the
characterisation of isocohomological homomorphisms.  We give several
equivalent characterisations.  The most attractive one may be that~$f$ is
isocohomological if and only if $f^*\colon \Der(B)\to\Der(A)$ is fully
faithful.  The general machinery developed above makes it easy to prove the
equivalence of these conditions.

Then we consider some classes of examples of isocohomological homomorphisms.
First we consider group convolution algebras, then crossed products.  Our
results on crossed products contain results on noncommutative tori as special
cases.  In the group case, we first exhibit how our definitions are related
to~\cite{Meyer:Combable_poly}.  We explain how to treat weighted variants of
the space $\Sch(G)$ for a discrete group~$G$.  I needed this generalisation in
a rather simple special case in~\cite{Meyer:Primes_Rep}.  Then we show that
the embedding $\C[G]\to\Sch(G)$ is isocohomological for discrete groups of
polynomial growth.  The proof reduces this to the special case where~$G$ is a
cocompact lattice in a connected nilpotent Lie group~$\bar{G}$.  We use the de
Rham complex of compactly supported differential forms on~$\bar{G}$ as a free
$\C[G]$\brd{}module resolution of the trivial representation of~$G$.  The same
argument shows that the embedding $\CCINF(\bar{G})\to\Sch(\bar{G})$ is
isocohomological.  For a Lie group~$G$, let $\Hol(G)$ be the algebra of
functions on~$G$ whose derivatives decay faster than any exponential in the
length function.  Our argument also shows that the embedding
$\CCINF(G)\to\Hol(G)$ is isocohomological for a large class of~$G$, including
all semi-simple Lie groups.  It seems likely that the Harish-Chandra-Schwartz
algebras of semi-simple Lie groups and their $p$\nbd{}adic analogues are also
isocohomological.  However, so far I could check this only for
$\mathrm{Sl}_2(\Q_p)$.

\section{Bornological modules over unital bornological algebras}
\label{sec:born_modules}

Throughout this article, we may consider bornological algebras over the real
or complex numbers.  We usually suppose that we work over~$\C$ in our
notation.  We write $X\inOb\Cat$ to denote that~$X$ is an object of a
category~$\Cat$.

\subsection{Bornologies}
\label{sec:bornologies}

We shall only meet complete convex bornological vector spaces in this article.
Therefore, we drop these qualifiers from our notation and tacitly assume all
bornologies to be complete and convex.  We refer
to~\cites{Hogbe-Nlend:Bornologies, Hogbe-Nlend:Completions,
  Prosmans-Schneiders:Born_Indlim, Meyer:Born_Top} for general results on
bornological vector spaces.

A \emph{bornology} on a vector space~$V$ is just a collection of
\emph{bounded} subsets.  Convexity of the bornology means that any bounded
subset is contained in an absolutely convex bounded subset.  If $T\subseteq V$
is an absolutely convex bounded subset, we let~$V_T$ be its linear span.
There is a unique seminorm on~$V_T$ whose unit ball is $T'\defeq
\bigcap_{\epsilon>0} (1+\epsilon)T$.  We call~$T$ a \emph{complete disk}
in~$V$ if $T=T'$ and~$V_T$ is a Banach space.  Completeness of the bornology
means that any bounded subset is contained in a complete disk.  The complete
disks in~$V$ form a directed set and $T\mapsto V_T$ is an inductive system of
Banach spaces with the additional property that the maps $V_S\to V_T$ for
$S\le T$ are injective.  We have a natural isomorphism of bornological vector
spaces $V\cong \varinjlim V_T$.  A closer look at this construction reveals
that the category of bornological vector spaces is equivalent to the category
of inductive systems of Banach spaces with injective structure maps
(see~\cite{Prosmans-Schneiders:Born_Indlim}).  Analysis in bornological vector
spaces is done by reduction to the Banach space case.  For instance, a
sequence converges in~$V$ if and only if it converges in~$V_T$ for some
complete disk~$T$.

We shall need the following examples of bornologies.  Any vector space~$V$ can
be written as a direct union of its finite dimensional subspaces.  This gives
rise to the \emph{fine bornology} on~$V$.  The precompact subsets of a Fréchet
space~$V$ form a bornology, called the \emph{precompact bornology} on~$V$.  We
remark that ``precompact'' is synonymous to ``totally bounded''; it is
equivalent to ``relatively compact'' in complete spaces.

Let $V$ and~$W$ be two bornological vector spaces.  We write $\Hom(V,W)$ for
the space of bounded linear maps $V\to W$.  A subset $S\subseteq \Hom(V,W)$ is
called \emph{equibounded} if $S(T)\subseteq W$ is bounded for any bounded
subset $T\subseteq V$.  The equibounded subsets form a complete convex
bornology on $\Hom(V,W)$.  The bornological vector space
$\Hom^{(n)}(V_1×\dotsb×V_n,W)$ of bounded $n$\nbd{}linear maps
$V_1×\dotsb×V_n\to W$ is defined similarly.  It is easy to see that there is a
bornological isomorphism
\begin{equation}  \label{eq:Hom_sums_products}
  \Hom\Bigl(\bigoplus_i X_i,\prod_j Y_j\Bigr) \cong \prod_{i,j} \Hom(X_i,Y_j).
\end{equation}
Moreover, composition is a bounded bilinear map
$$
\circ\colon \Hom(X,W)×\Hom(V,X)\to\Hom(V,W),\qquad (f,g)\mapsto f\circ g.
$$
Thus $\End(V)\defeq \Hom(V,V)$ is a bornological algebra.

The \emph{projective (complete bornological) tensor product} $V\hot W$ is a
(complete convex) bornological vector space with the universal property that
bounded bilinear maps $V\times W\to X$ into complete convex bornological
vector spaces correspond to bounded linear maps $V\hot W\to X$
(see~\cite{Hogbe-Nlend:Completions}).  There is a natural bornological
isomorphism
$$
\Hom(V_1\hot\dotsb\hot V_n,W) \cong \Hom^{(n)}(V_1\times\dotsb\times V_n,W).
$$
Hence~$\hot$ is commutative, associative, and satisfies $\C\hot
V\cong V\hot\C \cong V$ for all~$V$.  There are natural \emph{adjoint
  associativity} isomorphisms
\begin{equation}  \label{eq:primitive_adjoint_associativity}
  \Hom(V\hot W,X) \cong \Hom(V,\Hom(W,X)) \cong \Hom(W,\Hom(V,X))
\end{equation}
because these spaces are all bornologically isomorphic to $\Hom^{(2)}(V\times
W,X)$.  Equation~\eqref{eq:primitive_adjoint_associativity} implies
that~$\hot$ commutes with direct limits.  Especially,
\begin{equation}  \label{eq:hot_sums_products}
  \Bigl(\bigoplus_i X_i\Bigr)\hot \Bigl(\bigoplus_j Y_j\Bigr)
  \cong \bigoplus_{i,j} X_i\hot Y_j.
\end{equation}

In categories of topological vector spaces, adjoint associativity fails.
There are many useful topologies on spaces of continuous linear maps, but they
do not have particularly good algebraic properties.  In fact, the most
canonical structure on a space of continuous linear maps is a bornology,
namely, the equicontinuous bornology.  The complete projective topological
tensor product (\cite{Grothendieck:Produits}) cannot have a right or left
adjoint functor because it does not commute with direct sums and it does not
preserve kernels.  Instead, it commutes with direct products and preserves
cokernels.

\subsection{Bornological algebras and modules}
\label{sec:born_alg_mod}

A \emph{bornological algebra} is a bornological vector space together with a
bounded associative multiplication $A\times A\to A$.  A \emph{bornological
  (unital) left $A$\nbd{}module} is a bornological vector space~$M$ together
with a bounded bilinear multiplication map $A\times M\to M$ satisfying the
usual rules $a_1\cdot (a_2\cdot m)=(a_1\cdot a_2)\cdot m$, $1\cdot m=m$.  By
the universal property of the complete bornological tensor product, this data
is equivalent to a bounded linear map $A\hot M\to M$ satisfying similar
properties.  By adjoint
associativity~\eqref{eq:primitive_adjoint_associativity}, this is further
equivalent to a bounded unital algebra homomorphism $A\to\End(M)$.  The latter
description is, of course, not available for topological modules.  We let
$\Mod(A)$ be the \emph{category of bornological left $A$\nbd{}modules}, whose
objects are the bornological left $A$\nbd{}modules and whose morphisms are the
bounded $A$\nbd{}module homomorphisms.  In this section, we only consider
unital algebras and modules.  The important generalisation to essential
modules over quasi-unital algebras is treated in
Section~\ref{sec:quasi-unital_essential}.

Right modules are defined similarly using a homomorphism $A^\op\to\End(M)$, a
bilinear map $M\times A\to M$, or a bilinear map $M\hot A\to M$.  Here~$A^\op$
is the opposite algebra of~$A$.  Thus we write $\Mod(A^\op)$ for the category
of right $A$\nbd{}modules.  The category of $A,B$\brd{}bimodules is denoted by
$\Mod(A\hot B^\op)$.  We write $\Born$ for the category of bornological vector
spaces.  Thus $\Born=\Mod(\C)$.

Any morphism $f\colon M\to N$ in $\Mod(A)$ has both a kernel and a cokernel.
The kernel is simply the usual vector space kernel equipped with the subspace
bornology.  The cokernel is the quotient $N/\cl{f(M)}$, where the closure of
$f(M)$ is the smallest subspace of~$N$ that contains $f(M)$ and has the
property that any sequence in $\cl{f(M)}$ that converges in~$N$ has its limit
point again in $\cl{f(M)}$.  It is an important fact that this quotient is
again a complete convex bornological vector space
(see~\cite{Hogbe-Nlend:Bornologies}).  We also have canonical bornologies on
direct sums and products.  These fulfil the usual categorical requirements for
coproducts and products.  As a result, the category $\Mod(A)$ is complete and
cocomplete, that is, any diagram in $\Mod(A)$ has both a direct and an inverse
limit.  It is also clear that $\Mod(A)$ is an additive category.

\subsection{Adjoint associativity}
\label{sec:adjoint_associativity}

If $M,N\inOb\Mod(A)$, we let $\Hom_A(M,N)$ be the space of bounded
$A$\nbd{}module homomorphisms with the equibounded bornology.  Thus
$\Hom_A(M,N)$ is a bornological subspace of $\Hom(M,N)$.  If
$M\inOb\Mod(A^\op)$, $N\inOb\Mod(A)$, we define the \emph{$A$\nbd{}balanced
  projective tensor product} $M\hot_A N$ as the cokernel of the bounded linear
map
$$
M \hot A \hot N \to M\hot N,
\qquad
m\otimes a\otimes n\mapsto m\cdot a\otimes n- m\otimes a\cdot n.
$$
That is, we divide by the closure of the range of this map.  Thus $M\hot_A
N$ is a complete convex bornological vector space.

Now let $A,B,C$ be unital bornological algebras.  If $M\inOb\Mod(A\hot
B^\op)$ and $N\inOb\Mod(A\hot C^\op)$, then $\Hom_A(M,N)$ carries a
canonical $B,C$\brd{}bimodule structure by $b\cdot f\cdot c(m)\defeq f(m\cdot
b)\cdot c$ for all $b\in B$, $c\in C$, $m\in M$, $f\in\Hom_A(M,N)$.  Thus we
get a bifunctor
$$
\Hom_A \colon \Mod(A\hot B^\op)^\op \times \Mod(A\hot C^\op)
\to \Mod(B\hot C^\op).
$$
If $M\inOb\Mod(B\hot A^\op)$ and $N\inOb\Mod(A\hot C^\op)$, then
$M\hot_A N$ carries a canonical $B,C$\brd{}bimodule structure by $b\cdot
(m\otimes n)\cdot c\defeq (b\cdot m)\otimes (n\cdot c)$ for all $b\in B$,
$c\in C$, $m\in M$, $n\in N$.  Thus we get a bifunctor
$$
\hot_A \colon \Mod(B\hot A^\op) \times \Mod(A\hot C^\op)
\to \Mod(B\hot C^\op).
$$
Adjoint associativity~\eqref{eq:adjoint_associativity} generalises to natural
bornological isomorphisms
\begin{equation}  \label{eq:adjoint_associativity}
  \begin{aligned}
    \Hom_{B,C}(M\hot_A N,X) &\cong \Hom_{A,C}(N,\Hom_B(M,X)) \\
    &\cong \Hom_{B,A}(M,\Hom_C(N,X))
  \end{aligned}
\end{equation}
for $M\inOb\Mod(B\hot A^\op)$, $N\inOb\Mod(A\hot C^\op)$, $X\inOb\Mod(B\hot
C^\op)$.  The proof of~\eqref{eq:adjoint_associativity} identifies all three
spaces with the space of bounded bilinear maps $f\colon M×N\to X$ that satisfy
$f(b\cdot m\cdot a,n\cdot c)=b\cdot f(m,a\cdot n)\cdot c$ for all $b\in B$,
$m\in M$, $a\in A$, $n\in N$, $c\in C$.  Equations
\eqref{eq:Hom_sums_products} and~\eqref{eq:hot_sums_products} imply
\begin{equation}  \label{eq:HomA_hotA_sums_products}
  \begin{aligned}
    \Hom_A\Bigl(\bigoplus_i M_i,\prod_j N_j\Bigr) &\cong \prod_{i,j} \Hom_A(M_i,N_j),
    \\
    \Bigl(\bigoplus_i M_i\Bigr)\hot_A \Bigl(\bigoplus_j N_j\Bigr)
    &\cong \bigoplus_{i,j} M_i\hot_A N_j.
  \end{aligned}
\end{equation}

\subsection{Exact category structure}
\label{sec:ModA_exact}

In order to do homological algebra we need extensions.  There is a maximal
class of possible extensions in any additive category: a diagram
$K\overset{i}\to E\overset{p}\to Q$ with $p\circ i=0$ is called an
\emph{extension} if~$i$ is a kernel of~$p$ and~$p$ is a cokernel of~$i$.  In
our case, this simply means that the maps $K\to i(K)$ and $E/i(K)\to Q$ are
bornological isomorphisms with respect to the subspace bornology on $i(K)$ and
the quotient bornology on $E/i(K)$.  The notion of an \emph{exact category}
formalises the properties that a class of extensions should have.  One can
show that $\Mod(A)$ with the class of all extensions as defined above is an
exact category.  However, in most applications and, in particular, in this
article, we use a much smaller class of extensions.  An extension is called
\emph{linearly split} if there exists a bounded linear map $s\colon Q\to E$
such that $ps=\ID_Q$.  Equivalently, the forgetful functor $\Mod(A)\to\Born$
maps $K\into E\prto Q$ to a direct sum extension in $\Born$.  We always choose
this class of extensions in the following.

\begin{proposition}  \label{pro:ModA_exact}
  The additive category $\Mod(A)$ with the class of linearly split extensions
  is an exact category in the sense of Daniel Quillen.  Moreover, $\Mod(A)$ is
  complete and cocomplete.
\end{proposition}

\begin{proof}
  It is straightforward to verify that $\Mod(A)$ satisfies the axioms for an
  exact category in \cites{Keller:Appendix, Keller:Handbook}.  Completeness
  and cocompleteness mean that there are arbitrary direct and inverse limits.
  This follows from the existence of kernels, cokernels, direct sums, and
  direct products.
\end{proof}

Although trivial, Proposition~\ref{pro:ModA_exact} is important because
many definitions and theorems that are familiar from Abelian categories still
work in complete and cocomplete exact categories.

We use the following notation from~\cite{Keller:Appendix}.  The special
extensions $K\overset{i}\to E\overset{p}\to Q$ that are part of the exact
category structure are called \emph{conflations}; the maps~$i$ in conflations
are called \emph{inflations}, the maps~$p$ in conflations are called
\emph{deflations}.  In $\Mod(A)$, conflations are linearly split extensions,
inflations are linearly split bornological embeddings, and deflations are
linearly split bornological quotient maps.

Choosing a smaller class of extensions corresponds to doing \emph{relative}
homological algebra.  In our case, we work relative to the forgetful functor
$\Mod(A)\to\Born$, that is, we do homological algebra relative to the pure
analysis that occurs in the category $\Born$.  This explains why the bornology
creates no problems in the following algebraic constructions.

\subsection{Free and cofree modules}
\label{sec:free_cofree}

There is a natural isomorphism
\begin{equation}  \label{eq:HomAAM}
  \Hom_A(A,M)\cong M,  \qquad f\mapsto f(1),
\end{equation}
for any left or right bornological $A$\nbd{}module~$M$.  Together
with~\eqref{eq:adjoint_associativity} this yields natural bornological
isomorphisms
\begin{align}  \label{eq:AhotAM}
  A\hot_A M &\cong M,
  \\  \label{eq:NAhotA}
  N\hot_A A &\cong N,
  \\  \label{eq:HomAAhotVM}
  \Hom_A(A\hot V,M) &\cong \Hom(V,M),
  \\  \label{eq:HomAMHomAV}
  \Hom_A(M,\Hom(A,V)) \cong \Hom(M,V)
\end{align}
for any $M\inOb\Mod(A)$, $N\inOb\Mod(A^\op)$, $V\inOb\Born$.  Of
course, it is easy enough to prove \eqref{eq:AhotAM}--\eqref{eq:HomAMHomAV}
directly.  We call modules of the form $A\hot V$ \emph{free}, those of the
form $\Hom(A,V)$ \emph{cofree}.  Equations
\eqref{eq:HomAAhotVM}--\eqref{eq:HomAMHomAV} mean that the constructions of
free and cofree modules are left and right adjoints of the forgetful functor
$\Mod(A)\to\Born$.

\begin{lemma}  \label{lem:free}  \label{lem:cofree}
  Let $K\into E\prto Q$ be a conflation in $\Mod(A)$, let~$F$ be a free
  module and~$C$ a cofree module.  Then
  \begin{gather*}
    0\to\Hom_A(F,K)\to \Hom_A(F,E)\to \Hom_A(F,Q)\to 0,
    \\
    0\to \Hom_A(Q,C)\to\Hom_A(E,C)\to \Hom_A(K,C)\to 0
  \end{gather*}
  are conflations in $\Born$ and \emph{a fortiori} exact sequences of vector
  spaces.
\end{lemma}

\begin{proof}
  The assertion follows immediately from \eqref{eq:HomAAhotVM}
  and~\eqref{eq:HomAMHomAV}.
\end{proof}

An object~$P$ of an exact category~$\Cat$ is called \emph{projective} if the
covariant functor $\Hom_\Cat(P,\blank)$ is exact on conflations,
\emph{injective} if the contravariant functor $\Hom_\Cat(\blank,I)$ is exact
on conflations.  Thus Lemma~\ref{lem:free} implies that free modules are
projective and cofree modules are injective in $\Mod(A)$.  We say that an
exact category \emph{has enough projectives} if any object~$X$ admits a
deflation $P\to X$ with projective~$P$.  Dually, it \emph{has enough
  injectives} if any object~$X$ admits an inflation $X\to I$ with
injective~$I$.

\begin{proposition}
  \label{pro:ModA_enough_projectives}
  \label{pro:ModA_enough_injectives}
  The exact category $\Mod(A)$ has enough projectives and injectives.
\end{proposition}

\begin{proof}
  The map $A\hot M\overset{p}\to M$, $a\otimes m\mapsto a\cdot m$, is a
  deflation because it has the bounded linear section $m\mapsto 1\otimes m$.
  Its source is free and hence projective by Lemma~\ref{lem:free}.  Thus
  $\Mod(A)$ has enough projectives.  The map $M\overset{i}\to \Hom(A,M)$,
  $i(m)(a)\defeq a\cdot m$, is an inflation because it has the bounded linear
  section $f\mapsto f(1)$.  Its range is cofree and hence injective by
  Lemma~\ref{lem:cofree}.
\end{proof}

\section{Derived categories of bornological modules}
\label{sec:DerA}

\subsection{Chain complexes of bornological modules}
\label{sec:chains}

A chain complex over $\Mod(A)$ is given by $M = (M_m,\delta^M_m)_{m\in\Z}$,
where the~$M_m$ are bornological left $A$\nbd{}modules and the boundary maps
$\delta^M_m\colon M_m\to M_{m-1}$ are bounded left $A$\nbd{}module
homomorphisms satisfying $\delta^M_m\circ\delta^M_{m+1}=0$ for all $m\in\Z$.
We do not require chain complexes to be bounded below or above.  This has the
advantage that the category of chain complexes is still complete and
cocomplete.

We can also extend the functors $\Hom_A$ and $\hot_A$ to the level of chain
complexes.  Let $M = (M_m,\delta^M_m)_{m\in\Z}$ and $N =
(N_m,\delta^N_m))_{m\in\Z}$ be two chain complexes over $\Mod(A)$.  Let
$$
\Hom_A(M,N)_n\defeq \prod_{j\in\Z} \Hom_A(M_j,N_{n+j}).
$$
Together with the boundary maps $\delta(f)=\delta^N\circ f+ (-1)^{\abs{f}}
f\circ \delta^M$, where~$\abs{f}$ is the degree of~$f$, this defines a chain
complex of bornological vector spaces.
Let $M$ and~$N$ be chain complexes over $\Mod(A^\op)$ and $\Mod(A)$,
respectively.  Then
$$
(M\hot_A N)_n\defeq \bigoplus_{p+q=n} M_p\hot_A N_q
$$
with boundary map $\delta = \delta^M \hot_A\ID ± \ID\hot_A \delta^N$ with
appropriate signs defines a chain complex $M\hot_A N$ over $\Born$.  If $M$
and~$N$ are chain complexes of bornological $A,B$\brd{}bimodules and
$A,C$\brd{}bimodules, respectively, then $\Hom_A(M,N)$ is a chain complex of
bornological $B,C$\brd{}bimodules in a canonical way.  If $M$ and~$N$ are
chain complexes of bornological $B,A$\brd{}bimodules and $A,C$\brd{}bimodules,
respectively, then $M\hot_A N$ becomes a chain complex of bornological
$B,C$\brd{}bimodules.  There are natural isomorphisms of bornological chain
complexes
\begin{equation}  \label{eq:adjoint_associativity_chain}
  \begin{aligned}
    \Hom_{B,C}(M\hot_A N,X) &\cong \Hom_{A,C}(N,\Hom_B(M,X))
    \\ &\cong \Hom_{B,A}(M,\Hom_C(N,X))
  \end{aligned}
\end{equation}
if $M$, $N$, and~$X$ are chain complexes over $\Mod(B\hot A^\op)$, $\Mod(A\hot
C^\op)$, and $\Mod(B\hot C^\op)$, respectively.  This follows easily from the
corresponding assertion for modules in~\eqref{eq:adjoint_associativity} and
\eqref{eq:HomA_hotA_sums_products}.  Of course,
\eqref{eq:HomA_hotA_sums_products} remains valid for chain complexes.

\subsection{The homotopy category of chain complexes}
\label{sec:ho_chains}

We let $\Ho(A)$ be the category whose objects are the chain complexes over
$\Mod(A)$ and whose morphisms are given by $\Ho_A(M,N) \defeq
H_0(\Hom_A(M,N))$.  That is, $\Ho_A(M,N)$ is the space of homotopy classes of
chain maps, where chain maps and homotopies are required to be bounded and
$A$\nbd{}linear.  This is the \emph{homotopy category of chain complexes} over
$\Mod(A)$.  The homotopy category $\Ho(\Cat)$ can be constructed over any
additive category~$\Cat$.  It is well-known that $\Ho(\Cat)$ is a triangulated
category if~$\Cat$ is Abelian.  In fact, the axioms of a triangulated category
were introduced by Jean-Louis Verdier with the purpose of formalising some
properties of $\Ho(\Cat)$ that are needed to construct derived categories.  We
briefly recall the definition of the translation automorphism and the class of
exact triangles in $\Ho(\Cat)$.

The \emph{translation automorphism} $\Sigma\colon \Ho(\Cat)\to\Ho(\Cat)$ is
defined by $(\Sigma M)_m\defeq M_{m-1}$, $\delta^{\Sigma
  M}_m=-\delta^M_{m-1}$, and $(\Sigma f)_m = f_{m-1}$ for a chain map $f\colon
M\to N$.  The \emph{mapping cone} $\cone(f)$ of a chain map $f\colon M\to N$
is defined by $\cone(f)_m\defeq (N\oplus \Sigma M)_m$ with boundary map
$$
\delta^{\cone(f)} =
\begin{pmatrix}
  \delta^N & f \\ 0 & \delta^{\Sigma M}
\end{pmatrix}
$$
By construction, it fits into an extension of chain complexes $N\to
\cone(f) \to \Sigma M$, which splits if we disregard the boundary map.  The
resulting sequence of maps
$$
M\overset{f}\to N\to \cone(f) \to \Sigma M
$$
is called a \emph{mapping cone triangle}.  A morphism of triangles is a triple
of maps $(\mu,\nu,\gamma)$ giving rise to a commuting diagram
$$
\xymatrix{
  M \ar[r] \ar[d]^\mu & N\ar[r] \ar[d]^\nu & C \ar[r] \ar[d]^\gamma &
  \Sigma M \ar[d]^{\Sigma\mu}
  \\
  M' \ar[r] & N' \ar[r] & C' \ar[r] & \Sigma M'.
}
$$
This morphism is an isomorphism if $\mu$, $\nu$, and~$\gamma$ are
invertible.  A triangle in $\Ho(\Cat)$ is called \emph{exact} if it is
isomorphic to a mapping cone triangle.

A triangulated category is an additive category together with a translation
automorphism and a class of exact triangles satisfying some axioms
(see~\cites{Verdier:Thesis, Neeman:Triangulated}).  The axioms of a
triangulated category were invented to hold if~$\Cat$ is an Abelian category.
Since all the definitions above only need~$\Cat$ to be additive, it is not
surprising that the axioms still hold if~$\Cat$ is just an additive category.
This is routine to verify.  A detailed account can be found
in~\cite{Prosmans:Thesis}.  Readers with some background in operator algebra
K\nbd{}theory may also profit from reading the relevant sections
in~\cite{Meyer-Nest}.

\begin{proposition}  \label{pro:homotopy_category_triangulated}
  The category $\Ho(A)$ is a triangulated category in which every set of
  objects has a direct sum and a direct product.
\end{proposition}

\begin{proof}
  We have explained above why $\Ho(A)$ is a triangulated category.  Direct
  sums and products are easy to get: simply construct them in each degree
  separately and check that this has the required universal property.
\end{proof}

Let $F\colon \Mod(A)\to\Mod(B)$ be an additive functor, where~$B$ is another
bornological unital algebra.  Applying~$F$ to each entry of a chain complex,
we get a functor $F\colon \Ho(A)\to\Ho(B)$.  It is triangulated because it
commutes with the translation automorphism and the mapping cone construction.

Let $\Ho\defeq \Ho(\Born)$ be the homotopy category of bornological chain
complexes.  Recall that ${}^\op$ denotes opposite categories and algebras.

\begin{lemma}  \label{lem:chain_Hom_tensor}
  The bifunctors
  \begin{alignat*}{2}
    \Hom_A &\colon \Ho(A\hot B^\op)^\op×\Ho(A\hot C^\op)\to\Ho(B\hot C^\op),
    \\
    \hot_A &\colon \Ho(B\hot A^\op)×\Ho(A\hot C^\op)\to \Ho(B\hot C^\op)
  \end{alignat*}
  are triangulated functors in both variables and satisfy
  $$
  \Hom_A(\bigoplus_i M_i,\prod_j N_j) \cong \prod_{i,j} \Hom_A(M_i,N_j),
  \quad
  \bigoplus_i M_i\hot_A \bigoplus_j N_j \cong \bigoplus_{i,j} M_i\hot_A N_j.
  $$
\end{lemma}

\begin{proof}
  It is easy to see that $\Hom_A(M,N)$ and $M\hot_A N$ are functorial for
  bounded $A$\nbd{}linear chain maps and homotopies in each variable.  If we
  fix one variable, then the resulting one variable functor $\Hom_A$ or
  $\hot_A$ commutes with the translation automorphisms and the mapping cone
  construction.  Therefore, it is triangulated.  The last assertion is
  obvious.
\end{proof}

The following result is inspired by ideas of Marcel Bökstedt and Amnon
Neeman~\cite{Boekstedt-Neeman}.  It allows us to reduce assertions about chain
complexes to assertions about modules.  As usual, we embed
$\Mod(A)\subseteq\Ho(A)$ by sending a module to a chain complex concentrated
in degree~$0$.  This is a fully faithful embedding.

\begin{proposition}  \label{pro:ModA_generates_HoA}
  Let $\Subtri\subseteq\Ho(A)$ be a triangulated subcategory that contains
  $\Mod(A)$.  If~$\Subtri$ is closed under direct sums or under direct
  products, then $\Subtri=\Ho(A)$.
\end{proposition}

We may say, therefore, that $\Ho(A)$ is generated and cogenerated by the
subcategory $\Mod(A)$.  More generally, if~$\Cat$ is any additive category, we
may embed $\Cat\to\Ho(\Cat)$ as above.  Let $\Subtri\subseteq\Ho(\Cat)$ be a
triangulated subcategory that contains~$\Cat$.  The following proof still
works in this more general situation.  We get $\Subtri=\Ho(\Cat)$ if~$\Subtri$
is closed under countable direct sums and~$\Cat$ has kernels and countable
direct sums or, dually, if~$\Subtri$ is closed under countable direct products
and~$\Cat$ has quotients and countable direct products.

\begin{proof}
  First we prove that~$\Subtri$ contains all bounded chain complexes.  We do
  this by induction on the length of the chain complex.  Chain complexes of
  length~$1$ are all obtained by applying a power of the translation
  automorphism to a complex concentrated in degree~$0$.  Hence~$\Subtri$
  contains all complexes of length~$1$.  Consider some chain complex~$M$ of
  length~$n$, say, $M_{m+n-1}\to\dotsb\to M_m$.  There is an exact triangle
  involving~$M$, the truncated complex $M_{m+n-1}\to\dotsb\to M_{m+1}$,
  and~$M_m$ viewed as a complex concentrated in degree~$m$.  The last two
  belong to~$\Subtri$ by the induction hypothesis.  Since~$\Subtri$ is
  triangulated, $M$ belongs to~$\Subtri$ as well.  Thus~$\Subtri$ contains all
  bounded chain complexes.
  
  It remains to extend this to unbounded chain complexes.  Here we need the
  hypothesis on direct sums or products.  Suppose first that~$\Subtri$ is
  closed under direct sums.  Let $M=(M_n,\delta_n)$ be a possibly unbounded
  chain complex in $\Ho(A)$.  Truncating it, we obtain an inductive system of
  subcomplexes $(\Fil_nM)$ of the form
  $$
  0 \to M_n \to M_{n-1}\to \cdots \to M_{-n+1} \to \ker \delta_{-n} \to 0.
  $$
  Evidently, the original complex~$M$ is the direct limit of these
  subcomplexes.  Even more, the map $\bigoplus \Fil_nM\to M$ has an evident
  bounded linear section.  This implies that~$M$ is also equal to the
  \emph{homotopy} direct limit of the inductive system $(\Fil_nM)$.  That is,
  it fits into an exact triangle $\Sigma M\to \bigoplus \Fil_nM \to \bigoplus
  \Fil_nM \to M$.  Since $\Fil_nM$ is bounded for all $n\in\N$, we already
  know $\Fil_nM\in\Subtri$.  If~$\Subtri$ is triangulated and closed under
  direct sums, we obtain $M\in\Subtri$ as well.
  
  A dual argument works if~$\Subtri$ is closed under products.  Now we
  consider a projective system of quotient complexes $(\Fil^nM)$ of~$M$ of the
  form
  $$
  0 \to \coker \delta_{n+1} \to M_{n-1}\to \cdots \to M_{-n} \to 0.
  $$
  The original complex~$M$ is the inverse limit of these quotient
  complexes.  Since the map $M\to \prod \Fil^nM$ has an evident bounded linear
  section, $M$ is also a homotopy inverse limit of $(\Fil^nM)$.  Now proceed
  as above.
\end{proof}

\subsection{Exact chain complexes}
\label{sec:exact_Ho}

Let~$\Cat$ be an exact category, for instance, $\Mod(A)$.  A chain complex
$(M,\delta)$ over~$\Cat$ is \emph{exact in degree~$m$} if $\ker
(\delta_m\colon M_m\to M_{m-1})$ exists and $\delta_{m+1}\colon M_{m+1}\to
\ker\delta_m$ is a deflation (see~\cite{Keller:Handbook}).  Equivalently,
$\ker \delta_m$ and $\ker \delta_{m+1}$ exist and $\ker \delta_{m+1}
\overset{\subseteq}\longrightarrow M_{m+1}
\overset{\delta_{m+1}}\longrightarrow \ker \delta_m$ is a conflation.  A chain
complex is called \emph{exact} if it is exact in degree~$m$ for all $m\in\Z$.
A morphism $f\colon M\to N$ is called a \emph{quasi-isomorphism} if its
mapping cone is exact in the above sense.  In the case of $\Ho(A)$, a chain
complex of bornological left $A$\nbd{}modules is \emph{exact} in degree~$m$ if
and only if there exists a bounded, linear, not necessarily $A$\nbd{}linear
map $\sigma\colon \ker\delta_m\to M_{m+1}$ with $\delta_{m+1}\sigma=\ID$.  It
is exact if and only if there is a bounded linear contracting homotopy.  A
chain map is a quasi-isomorphism if and only if the forgetful functor maps it
to an isomorphism in $\Ho$.

\begin{proposition}  \label{pro:exact_triangulated}
  The subcategory of exact chain complexes $\Exact(A)\subseteq\Ho(A)$ is a
  thick, triangulated subcategory that is closed under direct sums and
  products.
\end{proposition}

\begin{proof}
  This is proved in great generality in~\cite{Neeman:Derived_Exact}.  The only
  hypothesis needed for $\Exact(\Cat)\subseteq\Ho(\Cat)$ to be thick and
  triangulated is that idempotent morphisms in~$\Cat$ should split, that is,
  have a kernel and a range object.  The assertion about direct sums and
  products follows easily if the class of conflations is closed under direct
  sums and products.
\end{proof}

We call $M\inOb\Ho(A)$ \emph{projective} if $\Ho_A(M,N)\cong0$ for all exact
chain complexes~$N$ (see~\cite{Keller:Handbook}).  The projectives form a
triangulated subcategory of $\Ho(A)$ that is closed under direct sums.  A
\emph{projective approximation} of $M\inOb\Ho(A)$ is a quasi-isomorphism
$f\colon P\to M$ with projective~$P$.

\begin{proposition}  \label{pro:projective_approx}
  Let~$\Cat$ be an exact category with split idempotents.  Projective
  approximations in $\Ho(\Cat)$ are unique up to isomorphism if they exist.
  The subcategory $\Subtri\subseteq\Ho(\Cat)$ of objects that admit a
  projective approximation is triangulated, and the construction of projective
  approximations defines a triangulated functor $P\colon \Subtri\to\Ho(\Cat)$.
  The same assertions hold for injective approximations.
\end{proposition}

\begin{proof}
  This is well-known.  See \cite{Meyer-Nest}*{Proposition 2.4} for a quick
  proof from the axioms of a triangulated category.  There it is assumed that
  $\Subtri=\Ho(\Cat)$.  It is easy to see that we still get the assertions
  above without this requirement.
\end{proof}

\begin{proposition}  \label{pro:HoA_enough_projectives}
  The category $\Ho(A)$ has enough projectives and injectives.
  
  There are triangulated functors $P,I\colon \Ho(\Cat)\to\Ho(\Cat)$ and
  natural transformations $P(M)\to M\to I(M)$ such that $P(M)\to M$ is a
  projective approximation and $M\to I(M)$ is an injective approximation for
  each $M\inOb\Ho(A)$.  The projective approximation functor commutes with
  direct sums, the injective approximation functor commutes with direct
  products.
\end{proposition}

\begin{proof}
  We only prove the existence of projective approximations.  The injective
  case is similar.  Let $\Subtri\subseteq\Ho(A)$ be the subcategory of objects
  that have a projective approximation.  It is evidently closed under direct
  sums and under the translation automorphism, and it is triangulated by
  Proposition~\ref{pro:projective_approx}.
  Proposition~\ref{pro:ModA_generates_HoA} yields that $\Subtri=\Ho(A)$
  once~$\Subtri$ contains $\Mod(A)$.  Therefore, it suffices to construct
  projective approximations for objects of $\Mod(A)$.  Using that $\Mod(A)$
  has enough projectives (Proposition~\ref{pro:ModA_enough_projectives}), we
  construct a ``projective resolution'' $P_\bullet\to M$ for any
  $M\inOb\Mod(A)$.  Viewed as a chain map $P_\bullet\to M$, this is a
  quasi-isomorphism because it is a resolution.  It is easy to see that a
  bounded below chain complex of projective modules is projective.  It follows
  from Proposition~\ref{pro:projective_approx} that the construction is
  functorial and defines a triangulated functor.  It clearly commutes with
  direct sums.
\end{proof}

\subsection{The derived category}
\label{app:born_derived}

The derived category of $\Mod(A)$ is defined as the localisation of $\Ho(A)$
at the subcategory $\Exact(A)$ of exact chain complexes.  Since there are
enough projectives and injectives, this localisation is easy to describe
explicitly.  It has the same objects as $\Ho(A)$ and morphisms
\begin{displaymath}
  \Der_A(M,N)
  \cong \Ho_A(P(M),N)
  \cong \Ho_A(P(M),I(N))
  \cong \Ho_A(M,I(N)).
\end{displaymath}
Here $P$ and~$I$ are the projective and injective approximation functors.  The
projective and injective approximation functors on $\Ho(A)$ descend to
triangulated functors $P,I\colon \Der(A)\to\Ho(A)$.  Their ranges are the
subcategories of projectives and injectives, respectively.  These
subcategories of $\Ho(A)$ are equivalent to $\Der(A)$ as triangulated
categories (see \cite{Meyer-Nest}*{Proposition 2.4}).

\begin{proposition}  \label{pro:free_generate_DerA}
  Let $\Subtri\subseteq\Der(A)$ be a triangulated subcategory.  Then
  $\Subtri=\Ho(A)$ provided~$\Subtri$ contains all free modules and is closed
  under direct sums or~$\Subtri$ contains all cofree modules and is closed
  under direct products.
\end{proposition}

\begin{proof}
  Suppose that~$\Subtri$ contains free modules and is closed under direct
  sums.  The other case is proved dually.  As in the proof of
  Proposition~\ref{pro:ModA_generates_HoA}, one shows that~$\Subtri$ contains
  all bounded below chain complexes of free modules.  The freeness of the
  entries of the chain complex is not affected by our truncations because we
  only have to truncate above.  Hence~$\Subtri$ contains a free resolution for
  any $M\inOb\Mod(A)$.  Since this free resolution becomes isomorphic to~$M$
  in $\Der(A)$, we get $\Mod(A)\subseteq\Subtri$.  Now the assertion follows
  from Proposition~\ref{pro:ModA_generates_HoA}.
\end{proof}

Let $F\colon \Ho(A)\to\Ho(B)$ be a covariant or contravariant triangulated
functor.  We may get such functors by extending an additive functor
$\Mod(A)\to\Mod(B)$.  Applying~$F$ to projective or injective approximations,
we get the left and right \emph{(total) derived functors} of~$F$, which are
denoted $\Left F$ and $\Right F$.  By convention, $\Left F(M)\defeq F\circ
P(M)$ and $\Right F(M)\defeq F\circ I(M)$ if~$F$ is covariant, and $\Left
F(M)\defeq F\circ I(M)$ and $\Right F(M)\defeq F\circ P(M)$ if~$F$ is
contravariant.  The natural maps $P(M)\to M\to I(M)$ induce natural
transformations $\Left F(M)\to F(M)\to \Right F(M)$ in both cases.  By
construction, the functors $\Left F$ and $\Right F$ descend to triangulated
functors $\Left F,\Right F\colon \Der(A)\to\Ho(B)$.

\begin{proposition}  \label{pro:preserve_projective}
  Let $F\colon \Mod(A)\to\Mod(B)$ be a functor that commutes with direct sums
  and preserves projectives.  Then the extended functor $F\colon
  \Ho(A)\to\Ho(B)$ also preserves projectives.
  
  Let $F\colon \Mod(A)\to\Mod(B)$ be a functor that commutes with direct
  products and preserves injectives.  Then the extended functor $F\colon
  \Ho(A)\to\Ho(B)$ also preserves injectives.
\end{proposition}

\begin{proof}
  The subcategory of projectives $M\inOb\Ho(A)$ for which $F(M)$ is again
  projective is a triangulated subcategory closed under direct sums.  By
  hypothesis, it contains all projective modules.  As in the proof of
  Proposition~\ref{pro:free_generate_DerA}, this implies that it contains all
  projective chain complexes.  The proof for injectives is dual.
\end{proof}

\subsection{Derived adjoint associativity}
\label{sec:derive_adjoint_associativity}

We now derive the bifunctors $\Hom_A(M,N)$ and $M\hot_A N$ of
Lemma~\ref{lem:chain_Hom_tensor}.

\begin{lemma}  \label{lem:derive_Hom_tensor}
  For any $M,N\inOb\Ho(A)$, the natural maps
  $$
  \Hom_A(P(M),N) \to \Hom_A(P(M),I(N)) \leftarrow \Hom_A(M,I(N))
  $$
  are chain homotopy equivalences, that is, isomorphisms in $\Ho$.  Similarly,
  $$
  P(M) \hot_A N \cong P(M) \hot_A P(N) \cong M \hot_A P(N)
  $$
  in $\Ho$ for all $M\inOb\Ho(A^\op)$, $N\inOb\Ho(A)$.
\end{lemma}

\begin{proof}
  We only prove that the map $\Hom_A(P(M),N)\to\Hom_A(P(M),I(N))$ is a chain
  homotopy equivalence.  The other cases are similar.  Let~$F$ be a free
  bornological left $A$\nbd{}module and let $E\inOb\Ho(A)$ be exact.  Then
  the chain complex $\Hom_A(F,E)$ has a bounded contracting homotopy by
  Lemma~\ref{lem:free}.  An exact sequence argument shows that the map
  $\Hom_A(F,N)\to \Hom_A(F,I(N))$ is an isomorphism in $\Ho$ because the
  mapping cone of $N\to I(N)$ is exact and the functor $\Hom_A(F,\blank)$ is
  triangulated.  Thus the assertion holds if~$M$ is a free module.
  Let~$\Subtri$ be the class of objects~$M$ for which the assertion holds for
  all $N\inOb\Ho(A)$.  Since our assertion only depends on $P(M)$, we may
  view this as a subcategory of $\Der(A)$.  It is triangulated and closed
  under direct sums by the properties of $P$ and $\Hom_A(\blank,N)$.  Hence
  the assertion follows from Proposition~\ref{pro:free_generate_DerA} and the
  special case of free modules treated above.
\end{proof}

Therefore, we may define covariant bifunctors
\begin{align*}
  \Right\Hom_A &\colon \Der(A\hot B^\op)^\op×\Der(A\hot C^\op)
  \to\Ho(B\hot C^\op) \prto \Der(B\hot C^\op),
  \\
  \Lhot_A &\colon \Der(B\hot A^\op)×\Der(A\hot C^\op)\to\Ho(B\hot C^\op)
  \prto \Der(B\hot C^\op)
\end{align*}
by
\begin{align*}
  \Right\Hom_A(M,N) &\defeq \Hom_A(P(M),N) \cong \Hom_A(P(M),I(N))
  \cong\Hom_A(M,I(N)),
  \\
  M \Lhot_A N &\defeq P(M) \hot_A N \cong P(M) \hot_A P(N) \cong
  M \hot_A P(N).
\end{align*}
They are triangulated in each variable and still satisfy the assertions about
direct sums and products in Lemma~\ref{lem:chain_Hom_tensor}.  They also
satisfy adjoint associativity, that is, there are natural isomorphisms (in
$\Ho$)
\begin{equation}  \label{eq:derived_adjoint_associativity}
  \begin{aligned}
    \Right\Hom_{B,C}(M\Lhot_A N,X) &\cong \Right\Hom_{A,C}(N,\Right\Hom_B(M,X))
    \\ &\cong \Right\Hom_{B,A}(M,\Right\Hom_C(N,X))
\end{aligned}
\end{equation}
for $M\inOb\Der(B\hot A^\op)$, $N\inOb\Der(A\hot C^\op)$, $X\inOb\Der(B\hot
C^\op)$.  For the proof, rewrite these chain complexes as
$\Right\Hom_{B,C}(P(M)\hot_A P(N),X)$,
$$
\Hom_{A,C}(P(N),\Hom_B(P(M),X)),
\qquad
\Hom_{B,A}(P(M),\Hom_C(P(N),X))
$$
with projective approximations $P(M)\to M$ and $P(N)\to N$ in $\Ho(B\hot
A^\op)$ and $\Ho(A\hot C^\op)$, respectively.
Thus~\eqref{eq:derived_adjoint_associativity} follows
from~\eqref{eq:adjoint_associativity_chain} provided $P(M)\hot_A
P(N)\inOb\Ho(B\hot C^\op)$ is projective.  This is evident for free modules
and follows in general from a by now familiar abstract nonsense argument
involving Proposition~\ref{pro:free_generate_DerA}.

\subsection{Extension theory}
\label{sec:Ext}

Compose the embedding $\Mod(A)\to\Ho(A)$ with the canonical functor
$\Ho(A)\to\Der(A)$ to get a functor $\Mod(A)\to\Der(A)$.  This functor is
still fully faithful because $\Der_A(M,N)=\Hom_A(P(M),N)\cong\Hom_A(M,N)$
if~$N$ is a chain complex concentrated in degree~$0$.  The spaces
$$
\Ext^n_A(M,N) \defeq \Der_A(M,\Sigma^n N)
$$
are of special importance.  They clearly vanish if $n<0$.  These extension
groups retain all the well-known properties in the case of Abelian categories.
Namely, if $K\into E\prto Q$ is a conflation in $\Mod(A)$, then there are
associated long exact sequences in $\Ext^n_A$ in both variables.  There is an
associative cup product
\begin{multline*}
  \Ext^n_A(M,N)×\Ext^m_A(L,M)\congto
  \Der_A(\Sigma^m M,\Sigma^{n+m} N)×\Der_A(L,\Sigma^m M)
  \\ \to \Der_A(L,\Sigma^{m+n} N)
  \congto \Ext^{m+n}_A(L,N)
\end{multline*}
for all $L,M,N\inOb\Mod(A)$, $n,m\in\N$.

There is an alternative description of $\Ext^n_A(M,N)$ by equivalence classes
of $n$\nbd{}step conflations of~$M$ by~$N$ as in classical Yoneda theory.  In
particular, $\Ext^1_A(M,N)$ is the space of all isomorphism classes of
conflations $N\into E\prto M$, where isomorphism means that there is an
isomorphism $E\cong E'$ that induces the identity maps on $N$ and~$M$.  All
this holds in any exact category with split idempotents.  It is easy to prove
for $\Mod(A)$ because there are enough projectives and injectives.

\section{Quasi-unital algebras and essential modules}
\label{sec:quasi-unital_essential}

In this section we treat bornological algebras that do not necessarily possess
a unit element.  We let~$A^+$ be the bornological algebra that we get by
adjoining a new unit element to~$A$.  Thus $A^+\defeq A\oplus\C\cdot 1$ as a
bornological vector space.  If $A\to\End(M)$ defines a (possibly non-unital)
module over~$A$, then the unique unital extension $A^+\to\End(M)$ defines a
unital module over~$A^+$.  As a result, the category of non-unital left
bornological $A$\nbd{}modules is isomorphic to $\Mod(A^+)$.  If we define
$\Hom_A$ and $\hot_A$ in the obvious fashion, then
$\Hom_A(M,N)=\Hom_{A^+}(M,N)$ and $M\hot_A N=M\hot_{A^+} N$.  Hence all the
constructions of the previous sections apply to categories of non-unital
modules without any change.  The only point where we have to be careful is
that the free and cofree modules are now $A^+\hot X$ and $\Hom(A^+,X)$.  In
general, $A$ is not projective as a left or right $A$\nbd{}module.

It is often desirable to work with a subcategory of essential modules
$\Mod(A)\subseteq\Mod(A^+)$.  For instance, let~$G$ be a locally compact group
and let $\CCINF(G)$ be the convolution algebra of smooth, compactly supported
functions on~$G$.  The category $\Mod(\CCINF(G))$ is identified
in~\cite{Meyer:Smooth} with the category of smooth representations of~$G$ on
bornological vector spaces.  In this section, we describe in general under
what assumptions on~$A$ a category of essential modules can be defined, and we
discuss some constructions with essential modules.  The main ideas for this
are already contained in~\cite{Meyer:Smooth}.  Then we introduce morphisms
between quasi-unital algebras and study the functors on module categories and
derived categories they induce.

\subsection{Basic definitions}
\label{sec:quasi-unital}

\begin{definition}[\cite{Meyer:Smooth}]  \label{def:approximate_identity}
  A bornological algebra~$A$ has an \emph{approximate identity} if for any
  bounded subset $S\subseteq A$ there is a sequence $(u_n)$ in~$A$ such that
  the sequences $u_n\cdot x$ and $x\cdot u_n$ converge to~$x$ uniformly for
  $x\in S$.  This means that the sequences convergence uniformly in the usual
  sense in the Banach space~$A_T$ for some complete bounded disk $T\subseteq
  A$ depending only on~$S$.
\end{definition}

\begin{definition}  \label{def:quasi_unital}
  A bornological algebra~$A$ is called \emph{quasi-unital} if it has an
  approximate identity and if~$A$ is projective both as a left and a right
  $A$\nbd{}module.
\end{definition}

Projectivity means that there exist a bounded left $A$\nbd{}module
homomorphism $l\colon A\to A^+\hot A$ and a bounded right $A$\nbd{}module
homomorphism $r\colon A\to A\hot A^+$ that are sections for the multiplication
maps.  Since any element $a\in A$ can be written as $\lim u_n\cdot a$ for some
approximate unit~$(u_n)$, we have $l(a)=\lim u_n\cdot l(a)$.  Therefore, the
range of~$l$ is contained in $A\hot A\subseteq A^+\hot A$.  The same argument
shows that~$r$ maps into $A\hot A$.

\begin{lemma}[\cite{Meyer:Smooth}*{Lemma 4.4}]  \label{lem:smoothening_sub}
  Let~$A$ be a quasi-unital algebra and let~$M$ be a bornological left
  $A$\nbd{}module.  Then the natural map $A\hot_A M\to A^+\hot_A M\cong M$ is
  always injective.
\end{lemma}

The following proposition is related to \cite{Meyer:Smooth}*{Proposition 4.7}.

\begin{proposition}  \label{pro:essential_module}
  Let~$A$ be a quasi-unital algebra, let~$M$ be a bornological left
  $A$\nbd{}module, and define $\mu\colon A\hot M\to M$ by $\mu(a\otimes
  m)\defeq am$.  Then the following conditions are equivalent:
  \begin{enumerate}[\ref{pro:essential_module}.1.]
  \item $\mu\colon A\hot M\to M$ is a bornological quotient map;

  \item $\mu\colon A\hot M\to M$ has a bounded linear section;
    
  \item the map $\mu_*\colon A\hot_A M\to M$ induced by~$\mu$ is a
    bornological isomorphism.

  \end{enumerate}
\end{proposition}

\begin{proof}
  Since $\mu(a_1a_2\otimes m)=\mu(a_1\otimes a_2m)$, it induces a map
  $\mu_*\colon A\hot_A M\to M$, which is injective by
  Lemma~\ref{lem:smoothening_sub}.  Hence it is a bornological isomorphism if
  and only if it is a bornological quotient map.  This is equivalent to~$\mu$
  itself being a bornological quotient map because $A\hot M\to A\hot_A M$ is a
  bornological quotient map.  Thus \ref{pro:essential_module}.1 and
  \ref{pro:essential_module}.3 are equivalent.  Trivially,
  \ref{pro:essential_module}.2 implies \ref{pro:essential_module}.1.
  Conversely, the right module map $r\colon A\to A\hot A$ induces a map
  $$
  r'\colon A\hot_A M \overset{r\hot_A\ID}\longrightarrow A\hot A\hot_A M
  \overset{\ID_A\hot\mu_*}\longrightarrow A\hot M
  $$
  such that $\mu\circ r'=\mu_*$.  Hence \ref{pro:essential_module}.3
  implies \ref{pro:essential_module}.2.
\end{proof}

\begin{definition}  \label{def:essential_modules}
  We call a module~$M$ over a quasi-unital algebra \emph{essential} if it
  satisfies the equivalent conditions of
  Proposition~\ref{pro:essential_module}.  We let $\Mod(A)\subseteq \Mod(A^+)$
  be the full subcategory of essential modules.
\end{definition}

\begin{example}  \label{exa:unital}
  Unital algebras are quasi-unital.  The unit is an approximate identity, and
  we can take $r(a)=1\otimes a$ and $l(a)\defeq a\otimes 1$.  An
  $A$\nbd{}module is essential if and only if it is unital.  Thus our notation
  is compatible with our previous definitions for unital algebras.
\end{example}

\begin{example}  \label{exa:CCINF}
  Let~$M$ be a smooth manifold and let $\CCINF(M)$ be the algebra of smooth
  functions on~$M$ with compact support and with the pointwise multiplication.
  This algebra is quasi-unital.  There is a sequence $(\phi_n)_{n\in\N}$ in
  $\CCINF(M)$ with $\sum \phi_n^2=1$.  Then $u_n\defeq \sum_{m\le n} \phi_m^2$
  is an approximate identity in $\CCINF(M)$ and we can take
  $$
  r(a)\defeq \sum_{n\in\N} \phi_n\otimes \phi_n a,
  \qquad
  l(a)\defeq \sum_{n\in\N} a\phi_n\otimes \phi_n.
  $$
  A $\CCINF(M)$\brd{}module~$M$ is essential if and only if for any bounded
  subset $S\subseteq M$ there is $n\in\N$ with $u_n\cdot m=m$ for all $m\in
  S$.
\end{example}

\begin{example}  \label{exa:CCINFG}
  The convolution algebra $\CCINF(G)$ for a locally compact group~$G$ is
  quasi-unital.  The category $\Mod(\CCINF(G))$ is isomorphic to the category
  of smooth representations of~$G$ on bornological vector spaces
  (\cite{Meyer:Smooth}).
\end{example}

\subsection{Smoothening and roughening functors}
\label{sec:smooth_rough}

\begin{definition}  \label{def:smoothening}
  Let~$A$ be a quasi-unital algebra.  The \emph{smoothening functor} and the
  \emph{roughening functor} for~$A$ are defined by
  \begin{alignat*}{2}
    \Smooth=\Smooth_A &\colon \Mod(A^+)\to \Mod(A^+),
    &\qquad M&\mapsto A\hot_A M,
    \\
    \Rough=\Rough_A &\colon \Mod(A^+)\to \Mod(A^+),
    &\qquad M &\mapsto \Hom_A(A,M).
  \end{alignat*}
\end{definition}

\begin{theorem}  \label{the:smoothening}
  The functor $\Smooth$ is an idempotent functor from $\Mod(A^+)$ onto the
  full subcategory $\Mod(A)$.  We have a natural bornological isomorphism
  $$
  \Hom_A(M,\Smooth N)\cong\Hom_A(M,N)
  \qquad \forall M\inOb\Mod(A),\ N\inOb\Mod(A^+).
  $$
  As a consequence, $\Smooth\colon \Mod(A^+)\to\Mod(A)$ is right adjoint to
  the embedding functor $\Mod(A)\to \Mod(A^+)$.
  
  When viewed as functors on $\Mod(A^+)$, the functors $\Smooth$ and $\Rough$
  satisfy
  $$
  \Smooth^2\cong \Smooth, \quad \Rough^2 \cong \Rough, \quad
  \Smooth\circ\Rough \cong \Smooth, \quad \Rough\circ\Smooth \cong \Rough.
  $$
  There is a natural isomorphism
  $$
  \Hom_A(\Smooth(M),N)\cong \Hom_A(M,\Rough(N))
  $$
  for all $M,N\inOb\Mod(A^+)$, so that $\Smooth$ is left adjoint to
  $\Rough$.
\end{theorem}

\begin{proof}
  The multiplication map $A\hot A\to A$ has a bounded linear section
  because~$A$ is quasi-unital.  Hence $A\hot_A A\cong A$ by
  Proposition~\ref{pro:essential_module}.  Since the balanced tensor product
  is associative, this means that $\Smooth$ is idempotent.  By definition,
  $\Smooth(M)\cong M$ if and only if~$M$ is essential.  Thus $\Smooth$ is a
  projection from $\Mod(A^+)$ onto $\Mod(A)$.  The map $\Smooth(N)\to N$ is
  bounded and injective for all~$N$ by Lemma~\ref{lem:smoothening_sub}.
  Therefore, the induced map $\Hom_A(M,\Smooth(N))\to\Hom_A(M,N)$ is bounded
  and injective as well.  If~$M$ is essential, then the natural map
  $\Smooth(M)\to M$ is a bornological isomorphism.  Since smoothening is
  functorial, any map $M\to N$ restricts to a map
  $M\cong\Smooth(M)\to\Smooth(N)$.  This provides a bounded section for the
  bounded injective map above, so that $\Hom_A(M,\Smooth(N))\cong\Hom_A(M,N)$.
  
  Adjoint associativity~\eqref{eq:adjoint_associativity} implies
  $\Hom_A(\Smooth(M),N)\cong \Hom_A(M,\Rough(N))$ for all $M,N\inOb\Mod(A^+)$.
  We already observed above that $\Smooth^2\cong\Smooth$.  By adjointness,
  this implies $\Rough^2\cong\Rough$.  Adjointness assertions also yield
  \begin{multline*}
    \Hom_A(M,\Rough\Smooth N)
    \cong \Hom_A(\Smooth M,\Smooth N)
    \\ \cong \Hom_A(\Smooth M,N)
    \cong \Hom_A(M,\Rough N)
  \end{multline*}
  for all $M,N\inOb\Mod(A^+)$.  Hence $\Rough\Smooth\cong\Rough$.  Finally,
  If $M\inOb\Mod(A)$, $N\inOb\Mod(A^+)$, then
  \begin{multline*}
    \Hom_A(M,\Smooth \Rough N)
    \cong \Hom_A(M,\Rough N)
    \cong \Hom_A(\Smooth M, N)
    \\ \cong \Hom_A(M,N)
    \cong \Hom_A(M,\Smooth N).
  \end{multline*}
  Since both $\Smooth\Rough N$ and $\Smooth N$ belong to $\Mod(A)$, this
  implies $\Smooth\Rough\cong\Smooth$.
\end{proof}

\begin{corollary}  \label{cor:smooth_exact}
  The smoothening functor $\Smooth\colon \Mod(A^+)\to\Mod(A^+)$ is exact and
  preserves projectives, and it commutes with arbitrary direct limits in the
  category $\Mod(A^+)$.
\end{corollary}

\begin{proof}
  The functors $\Rough$ and $\Smooth$ are exact on $\Mod(A^+)$ because $A$ is
  projective as a right $A$\nbd{}module.  Since the right adjoint $\Rough$
  of~$\Smooth$ is exact, $\Smooth$ maps projectives in $\Mod(A^+)$ to
  projectives in $\Mod(A^+)$.  Furthermore, $\Smooth$ commutes with direct
  limits in $\Mod(A^+)$ because it has a right adjoint.
\end{proof}

\begin{theorem}  \label{the:essential_hereditary}  
  The subcategory $\Mod(A)\subseteq\Mod(A^+)$ is exact and hereditary for
  subspaces and quotients.  That is, if $K\into E\prto Q$ is a conflation in
  $\Mod(A^+)$, then~$E$ is essential if and only if $K$ and~$Q$ are both
  essential.
\end{theorem}

\begin{proof}
  Let $K\into E\prto Q$ be a conflation in $\Mod(A^+)$.  Since smoothening is
  exact, the rows in the commuting diagram
  $$
  \xymatrix{
    \Smooth K \ar[r] \ar[d] &
    \Smooth E \ar[r] \ar[d] &
    \Smooth Q \ar[d] \\
    K \ar[r] & E\ar[r] & Q
  }
  $$
  are conflations.  Suppose first that~$E$ is essential.  It is easy to see
  that quotients of essential modules are again essential.  Hence the maps
  $\Smooth E\to E$ and $\Smooth Q\to Q$ are isomorphisms by
  Theorem~\ref{the:smoothening}.  It follows that the map $\Smooth K\to K$ is
  an isomorphism as well, that is, $K$ is essential.  Thus $\Mod(A)$ is
  hereditary for subspaces and quotients.  Suppose conversely that $K$ and~$Q$
  are essential.  Then the maps $\Smooth K\to K$ and $\Smooth Q\to Q$ are
  isomorphisms.  This implies that the map $\Smooth E\to E$ is an isomorphism
  by the Five Lemma.  Thus~$E$ is essential.
\end{proof}

\subsection{Derived categories of essential modules}
\label{sec:derived_essential}

\begin{theorem}  \label{the:essential_ModA_properties}
  Let~$A$ be a quasi-unital bornological algebra.  The category $\Mod(A)$ of
  essential bornological left $A$\nbd{}modules with the class of linearly
  split extensions is a complete, cocomplete exact category.  It has enough
  projectives and injectives.
  
  The embedding $\Mod(A)\to\Mod(A^+)$ is exact, preserves projectives, and
  commutes with direct limits and finite inverse limits.
\end{theorem}

\begin{proof}
  The category $\Mod(A)$ is an exact subcategory of $\Mod(A^+)$ by
  Theorem~\ref{the:essential_hereditary} (that is, it is closed under
  extensions) and hence an exact category in its own right.  We already know
  that $\Mod(A^+)$ is complete and cocomplete.  The subcategory $\Mod(A)$ is
  evidently closed under direct sums in $\Mod(A^+)$.  It is also hereditary
  for subspaces and quotients by Theorem~\ref{the:essential_hereditary}.
  Therefore, $\Mod(A)$ is closed under arbitrary direct limits and finite
  inverse limits.  The inverse limit $\varprojlim M \in \Mod(A^+)$ of a
  diagram~$M$ in $\Mod(A)$ need not belong to $\Mod(A)$ any more.  However,
  $\Smooth \varprojlim M$ has the right universal property for an inverse
  limit in the subcategory $\Mod(A)$ by Theorem~\ref{the:smoothening}.  Hence
  $\Mod(A)$ is still a complete, cocomplete exact category and the embedding
  into $\Mod(A^+)$ commutes with direct limits and finite inverse limits, but
  not necessarily with infinite inverse limits.
  
  Since $\Mod(A^+)$ has enough projectives and injectives, any $M\inOb\Mod(A)$
  has a projective resolution $P\prto M$ and an injective coresolution $M\into
  I$ in $\Mod(A^+)$.  By Corollary~\ref{cor:smooth_exact}, $\Smooth P\prto M$
  is a projective resolution in $\Mod(A)$.  Even more, $\Smooth P$ is still
  projective in $\Mod(A^+)$.  Therefore, $\Mod(A)$ has enough projectives and
  the embedding into $\Mod(A^+)$ preserves projectives.  Since $\Smooth$ is
  right adjoint to an exact functor, it maps injectives in $\Mod(A^+)$ to
  injectives in $\Mod(A)$ (however, these need not be injective in $\Mod(A^+)$
  any more).  Since $\Smooth$ is also exact, $M\into \Smooth I$ is an
  injective coresolution of~$M$ in $\Mod(A)$.  Hence $\Mod(A)$ has enough
  injectives.
\end{proof} 

Let $R\subseteq \Mod(A^+)$ be the subcategory of rough modules.  Then
$\Smooth\colon R\to \Mod(A)$ and $\Rough\colon \Mod(A)\to R$ are isomorphisms
of exact categories.  Hence there is no need to study the category of rough
modules.

We let $\Ho(A)$ be the homotopy category of unbounded chain complexes over
$\Mod(A)$.  Let $\Exact(A)\subseteq \Ho(A)$ be the subcategory of exact chain
complexes and let $\Der(A)$ be the localisation of $\Ho(A)$ at $\Exact(A)$.
Theorem~\ref{the:essential_ModA_properties} allows to extend the assertions
about $\Ho(A)$ and $\Der(A)$ for unital~$A$ in Section~\ref{sec:DerA} to the
case of quasi-unital~$A$.  That is, $\Ho(A)$ is still a triangulated category
with direct sums and products, and $\Exact(A)$ is a localising subcategory.
The category $\Ho(A)$ has enough projectives and injectives, and projective
and injective approximations are functorial.  Analogues of Propositions
\ref{pro:ModA_generates_HoA} and~\ref{pro:free_generate_DerA} also hold, if we
replace free and cofree modules in the statement of
Proposition~\ref{pro:free_generate_DerA} by the special projective and
injective modules of the form $A\hot X$ and $\Smooth \Hom(A,X)$.

Since we can find projective resolutions in the category $\Mod(A^+)$ that lie
in the subcategory $\Mod(A)$, the natural map $\Der(A)\to\Der(A^+)$ is fully
faithful.  Since $\Smooth$ and $\Rough$ are exact functors, they give rise to
triangulated functors on $\Ho(A^+)$ and $\Der(A^+)$.  All the assertions of
Theorem~\ref{the:smoothening} remain valid in these categories.

We can modify the statement of adjoint associativity so that all involved
objects again belong to the appropriate derived categories of \emph{essential}
modules.  One checks easily that $M\hot_A N$ is an essential
$B,C$\brd{}bimodule if $M$ and~$N$ are essential $B,A$- and
$A,C$\brd{}bimodules.  This implies the corresponding assertion for chain
complexes and hence for $M\Lhot_A N$.  Using the properties of the
smoothening, we easily get
\begin{equation}  \label{eq:derived_adjoint_associativity_II}
  \begin{aligned}
    \Right\Hom_{B,C}(M\Lhot_A N,X)
    &\cong \Right\Hom_{A,C}(N,\Smooth_{A\hot C^\op} \Right\Hom_B(M,X))
    \\ &\cong \Right\Hom_{B,A}(M,\Smooth_{B\hot A^\op} \Right\Hom_C(N,X))
 \end{aligned}
\end{equation}
for all $M\inOb\Der(B\hot A^\op)$, $N\inOb\Der(A\hot C^\op)$, $X\inOb\Der(B\hot
C^\op)$.

\subsection{Multiplier algebras and Morita equivalence}
\label{sec:mult}

Let~$A$ be a quasi-unital bornological algebra.  The \emph{left and right
  multiplier algebras} $\Mult_l(A)$ and $\Mult_r(A)$ are defined as the
algebras of bounded right and left $A$\nbd{}module homomorphisms $A\to A$,
respectively.  The \emph{multiplier algebra} of~$A$ is the algebra of pairs
$(L,R)$ consisting of a left and a right multiplier, such that $a\cdot L(b) =
R(a)\cdot b$ for all $a,b\in A$.  It is a closed unital subalgebra of
$\Mult_l(A)\oplus \Mult_r(A)^\op$.  There is a canonical map $A\to\Mult(A)$,
which is injective if~$A$ is quasi-unital (see also~\cite{Meyer:Smooth}).

For example, if~$G$ is a Lie group, then the multiplier algebra of the
convolution algebra $\CCINF(G)$ is the algebra of compactly supported
distributions on~$G$ (\cite{Meyer:Smooth}).  The multiplier algebra of
$\CCINF(M)$ with the pointwise multiplication is $\CINF(M)$, the space of
smooth functions on~$M$ without support restriction.  If~$A$ is unital, then
$A\cong\Mult(A)$.

By definition, $A$ is a $\Mult_l(A),A$\nbd{}bimodule.  Hence $\Smooth(M)\defeq
A\hot_A M$ carries a canonical $\Mult_l(A)$\brd{}module structure for any
$M\inOb\Mod(A^+)$.  Thus essential $A$\nbd{}modules are at the same time also
bornological unital $\Mult_l(A)$\brd{}modules.  Similarly, essential right
$A$\nbd{}modules are bornological unital $\Mult_r(A)$\brd{}modules.  Both
kinds of modules are bornological unital modules over $\Mult(A)$ and
$\Mult(A)^\op$, respectively.  Thus we get a functor $E\colon
\Mod(A)\to\Mod(\Mult(A))$.  We also have a forgetful functor $F\colon
\Mod(\Mult(A))\to\Mod(A)$, $M\mapsto A\hot_A M$.  One can show that $A\hot_A
M\cong A\hot_{\Mult(A)} M$.

\begin{proposition}  \label{pro:fully_faithful_modules}
  The functor $E\colon \Mod(A)\to\Mod(\Mult(A))$ is left adjoint to the
  functor $F\colon \Mod(\Mult(A))\to\Mod(A)$.  The functor~$E$ is exact, fully
  faithful, commutes with direct limits, and preserves projectives.  The
  functor $E\colon \Der(A)\to\Der(\Mult(A))$ obtained from~$E$ is also fully
  faithful.
\end{proposition}

\begin{proof}
  It is clear that~$E$ is fully faithful.  The adjointness of $E$ and~$F$
  therefore follows from the universal property of the smoothening.  The
  adjointness of $E$ and~$F$ implies that~$E$ preserves projectives and
  commutes with direct limits because~$F$ is exact.  Exactness of~$E$ is
  trivial and implies that~$E$ gives rise to a functor
  $\Der(A)\to\Der(\Mult(A))$.  The latter is still fully faithful because~$E$
  is fully faithful and preserves projectives.
\end{proof}

The above construction is not an equivalence of categories.  For instance,
$\Mult(A)$ itself is not an essential $A$\nbd{}module.  We have $\Smooth
\Mult(A)=A$ because
$$
\Smooth \Mult(A)\subseteq \Smooth \Mult_l(A) = \Smooth \Rough(A) = A.
$$

One may wonder which category one should study, unital $\Mult(A)$\brd{}modules
or essential $A$\nbd{}modules.  I prefer the latter setting.  One reason for
this is that $\Mult(A)$ can become rather unwieldy.  This happens, for
instance, if $A=\CCINF(G)$ for a reductive group over the adele ring of~$\Q$.
Moreover, some important Morita equivalences work for $\Mod(A)$ but not for
$\Mod(\Mult(A))$.

Let $A$ and~$B$ be quasi-unital bornological algebras.  A \emph{Morita
  equivalence} between $A$ and~$B$ is given by essential bimodules
$P\inOb\Mod(A\hot B^\op)$, $Q\inOb\Mod(B\hot A^\op)$ together with bounded
bimodule isomorphisms $P\hot_B Q\cong A$, $Q\hot_A P\cong B$, such that~$P$ is
projective as a $B$\nbd{}module and~$Q$ is projective as an $A$\nbd{}module.

Given such bimodules, we define functors $M\mapsto P\hot_B M$ and $M\mapsto
Q\hot_A M$ between $\Mod(A)$ and $\Mod(B)$.  By hypothesis, these two functors
are inverse to each other.  Hence we get an equivalence of categories
$\Mod(A)\cong\Mod(B)$.  The projectivity hypothesis on $P$ and~$Q$ ensures
that these two functors are exact.  It is conceivable that the projectivity of
$P$ and~$Q$ holds automatically.  I have not seriously investigated Morita
equivalence and present it here only to exhibit a problem with the category
$\Mod(\Mult(A))$.

Consider the following rather simple example.  Let $A=\C$ and let $B$ the
Banach algebra $\ell_1(H)$ of trace class operators on a separable Hilbert
space~$H$.  Equip~$\ell_1(H)$ with the von Neumann bornology.  One can check
that the Hilbert space~$H$ viewed as column and row vectors gives a Morita
equivalence between $\C$ and $\ell_1(H)$ in the above sense.  Hence the
categories $\Mod(\C)$ and $\Mod(\ell_1(H))$ are equivalent.  The multiplier
algebra of~$\ell_1(H)$ is the algebra of all bounded operators on~$H$.  Its
canonical bornology is equal to the von Neumann bornology coming from the
usual operator norm.  The compact operators on~$H$ form a closed two-sided
ideal in $\Mult(\ell_1(H))$.  Since the ideal structure is an invariant under
Morita equivalence, $\Mult(\ell_1(H))$ cannot be Morita equivalent to~$\C$.
Hence the module categories $\Mod(\C)$ and $\Mod(\Mult(\ell_1(H)))$ are not
equivalent.

\subsection{Morphisms between quasi-unital bornological algebras}
\label{sec:morphisms}

Throughout this section, we let $A$ and~$B$ be quasi-unital bornological
algebras.  A bounded homomorphism $f\colon A\to B$ can always be extended to a
bounded unital homomorphism $f^+\colon A^+\to B^+$ and hence induces a functor
$f^*\defeq (f^+)^*\colon \Mod(B^+)\to\Mod(A^+)$.  However, this functor need
not map $\Mod(B)$ into $\Mod(A)$.  A necessary condition for this is that
$f^*(B)\in \Mod(A)$, that is, $A\hot_A B\cong B$.  In the following, we
frequently write~$B$ instead of $f^*(B)$ to avoid clutter.  We call a bounded
homomorphism $f\colon A\to B$ \emph{essential} if~$B$ is an essential
$A$\nbd{}bimodule, that is,
$$
A\hot_A B\hot_A A \cong A\hot_A B\cong B\hot_A B \cong B.
$$
More generally, we may consider a bounded algebra homomorphism $f\colon
A\to\Mult(B)$.  This still allows us to view~$B$ as an $A$\nbd{}bimodule.  We
call~$f$ a \emph{morphism} and write $f\colon A\dashrightarrow B$ if~$B$ is an
essential $A$\nbd{}bimodule.  A morphism $f\colon A\dashrightarrow B$ is
called \emph{proper} if~$f$ is a bounded map $A\to B$.  Notice that not all
bounded algebra homomorphisms $A\to B$ are proper morphisms in the above
sense: we require, in addition, that~$B$ should be an essential
$A$\nbd{}bimodule.

If~$A$ is unital, then a morphism $f\colon A\dashrightarrow B$ is the same as
a bounded unital algebra homomorphism $A\to\Mult(B)$.  Hence morphisms
$A\dashrightarrow B$ are the same as bounded unital algebra homomorphisms if
both $A$ and~$B$ are unital.

\begin{lemma}  \label{lem:morphism_Mult}
  Any morphism $f\colon A\dashrightarrow B$ extends uniquely to a bounded
  unital algebra homomorphism $\bar{f}\colon \Mult(A)\to\Mult(B)$.
  
  If $f\colon A\dashrightarrow B$ and $g\colon B\dashrightarrow C$ are
  morphisms, then $\bar{g}\circ f\colon A\to \Mult(B)\to\Mult(C)$ is a
  morphism $A\dashrightarrow C$, denoted $g\circ f$.  The quasi-unital
  bornological algebras form a category with this composition of morphisms.
\end{lemma}

\begin{proof}
  Let $f\colon A\to\Mult(B)$ be a morphism.  Since elements of the form
  $a\cdot b$ with $a\in A$, $b\in B$, span a dense subspace of~$B$, the
  extension~$\bar{f}$ is unique if it exists.  We can view~$B$ as a unital
  bimodule over $\Mult(A)$ because it is an essential $A$\nbd{}module.  This
  defines a bounded unital algebra homomorphism $\bar{f}\colon
  \Mult(A)\to\Mult(B)$ extending~$f$.
  
  To check that $\bar{g}\circ f$ is a morphism, we must show that~$C$ is an
  essential bimodule over~$A$.  Since $f$ and~$g$ are morphisms, we get
  \begin{multline*}
    A\hot_A C\hot_A A
    \cong A\hot_A (B\hot_B C\hot_B B)\hot_A A
    \\ \cong (A\hot_A B) \hot_B C \hot_B (B\hot_A A)
    \cong B\hot_B C\hot_B B \cong C
  \end{multline*}
  as desired.  It is easy to see that the composition of morphisms defined
  above gives rise to a category.
\end{proof}

We get a functor from the category of quasi-unital algebras to the category of
unital algebras by sending $A\mapsto \Mult(A)$, $f\mapsto\bar{f}$.  This
functor is faithful because the canonical map $A\to\Mult(A)$ is injective for
algebras with an approximate identity.  However, it is not fully faithful,
that is, there are bounded unital algebra homomorphisms $\Mult(A)\to\Mult(B)$
that do not come from a morphism $A\dashrightarrow B$.

The following lemma is useful to check that an algebra is quasi-unital.

\begin{lemma}  \label{lem:embed_quasiunital}
  Let $f\colon A\to B$ be a bounded homomorphism between two bornological
  algebras.  If~$A$ is quasi-unital and $f^*(B)$ is an essential
  $A$\nbd{}bimodule, then~$B$ is also quasi-unital and~$f$ is a proper
  morphism $A\to B$.
\end{lemma}

\begin{proof}
  Let $S\subseteq B$ be bounded.  There exist bounded disks $S_A\subseteq A$
  and $S_B\subseteq B$ such that $S\subseteq f(S_A)\cdot S_B\cdot f(S_A)$
  because the multiplication map $A\hot B\hot A\to B$ is supposed to be a
  bornological quotient map.  Since~$A$ is quasi-unital, we can find a
  sequence $(u_n)$ in~$A$ that acts as an approximate identity on~$S_A$.  Then
  $f(u_n)$ serves as an approximate identity on $f(S_A)\cdot S_B\cdot f(S_A)$
  and hence also on~$S$.  Thus~$B$ has an approximate identity.  The right
  $A$\nbd{}module section $r\colon A\to A\hot A$ for the multiplication map
  induces a bounded right $B$\nbd{}module homomorphism
  $$
  r_*\colon B\cong A\hot_A B
  \overset{r}\longrightarrow A\hot A\hot_A B
  \cong A\hot B
  \overset{f\hot\ID_B}\longrightarrow B\hot B
  $$
  with $\mu_B\circ r_*=\ID_B$.  A left module section is constructed
  similarly.
\end{proof}

\subsection{Functoriality of module categories and derived categories}
\label{sec:functoriality}

We let $f\colon A\dashrightarrow B$ be a morphism, described by a bounded
algebra homomorphism $f\colon A\to\Mult(B)$.  We may view~$B$ as a
$B,A$\brd{}bimodule or as an $A,B$\brd{}bimodule via~$f$.  This yields four
additive functors
\begin{alignat*}{2}
  f_!&\colon \Mod(A)\to\Mod(B),
  &\qquad f_!(M) &\defeq B\hot_A M,
  \\
  f_*&\colon \Mod(A)\to\Mod(B),
  &\qquad f_*(M)&\defeq \Smooth_B \Hom_A(B,M),
  \\
  f^*&\colon \Mod(B)\to\Mod(A),
  &\qquad f^*(M) &\defeq B\hot_B M,
  \\
  f^!&\colon \Mod(B)\to\Mod(A),
  &\qquad f^!(M)&\defeq \Smooth_A \Hom_B(B,M).
\end{alignat*}
By the definition of essential modules, $f^*(M)$ is equal to~$M$ viewed as a
left $A$\nbd{}module via $A\to\Mult(B)\to\End(M)$, where the second map is the
unique bounded extension of the $B$\nbd{}module structure $B\to\End(M)$.

\begin{example}
  Let~$G$ be a locally compact group and let $H\subseteq G$ be a closed
  subgroup.  Then we get an induced homomorphism $\CCINF(H)\to\CCINF(G)$ by
  sending $f\in\CCINF(H)$ to the multiplier of $\CCINF(G)$ associated to the
  compactly supported distribution $\psi\mapsto \int \psi(h)
  f(h)\,d\lambda_H(h)$ on~$G$, where $d\lambda_H$ denotes a left invariant
  Haar measure on~$H$.  This defines a morphism $i_H^G\colon \CCINF(H)
  \dashrightarrow \CCINF(G)$.  It is proper if and only if~$H$ is an open
  subgroup of~$G$.  The category $\Mod(\CCINF(G))$ is isomorphic to the
  category of smooth representations of~$G$ (\cite{Meyer:Smooth}).  Under this
  isomorphism, the functors $(i_H^G)^*$, $(i_H^G)_*$, and $(i_H^G)_!$
  correspond to the restriction functor, the induction functor, and the
  compact induction functor twisted by a relative modular function
  (see~\cite{Meyer:Smooth}).  The functor $(i_H^G)^!$ is called rough
  restriction functor in~\cite{Meyer:Smooth} and does not agree with any
  classical construction.
\end{example}

\begin{proposition}  \label{pro:adjointness_morphism}
  Let $f\colon A\dashrightarrow B$ be a morphism.  Then there are natural
  bornological isomorphisms
  $$
  \Hom_B(N,f_*(M))\cong\Hom_A(f^*(N),M),
  \qquad
  \Hom_B(f_!(M),N)\cong\Hom_A(M,f^!(N))
  $$
  for all $M\inOb\Mod(A)$, $N\inOb\Mod(B)$.  Thus~$f^*$ is left adjoint
  to~$f_*$ and~$f_!$ is left adjoint to~$f^!$.
\end{proposition}

\begin{proof}
  Use adjoint associativity and the universal property of $\Smooth$.
\end{proof}

The following proposition uses the functor $f^*\colon \Mod(B^+)\to\Mod(A^+)$
induced by a proper morphism $f\colon A\to B$.  Notice that its restriction to
the subcategory $\Mod(B)$ agrees with the functor~$f^*$ defined above.

\begin{proposition}  \label{pro:adjointness_proper_morphism}
  Let $f\colon A\to B$ be a proper morphism.  Then the natural maps
  \begin{gather*}
    \Smooth_A (f^* M) = A\hot_A M \to B\hot_B M = f^*(\Smooth_B M),
    \\
    f^*(\Rough_B M) = \Hom_B(B,M) \to \Hom_A(A,f^*M) = \Rough_A f^*(M),
  \end{gather*}
  induced by $f\colon A\to B$ are bornological isomorphisms for all
  $M\inOb\Mod(B^+)$.  Therefore, $f^!=f^*$.
\end{proposition}

\begin{proof}
  It follows from Lemma~\ref{lem:smoothening_sub} that the natural bounded
  $A$\nbd{}module homomorphisms $\Smooth_A f^*(M)\to f^*M$ and $f^*(\Smooth_B
  M)\to f^* M$ are injective.  We claim that $\Smooth_A (f^*M)$ and
  $f^*(\Smooth_B M)$ are the same subspace of $f^*M$ with the same bornology.
  We have a bounded inclusion $\Smooth_A (f^*M) = A\hot_A M \to B\hot_B M =
  f^*(\Smooth_B M)$.  Since $f^*(\Smooth_B M)$ is an essential
  $A$\nbd{}module, Theorem~\ref{the:smoothening} yields that the map
  $f^*(\Smooth_B M)\to f^*M$ factors through $\Smooth_A f^*(M)$.  Hence
  $\Smooth_A f^*(M)=f^*(\Smooth_B M)$ as asserted.
  
  We use this and the general properties of smoothening and roughening
  functors to prove $f^*(\Rough_B M) \cong \Rough_A f^*(M)$.  We have natural
  isomorphisms
  \begin{multline*}
    \Rough_A f^*(M)
    \cong \Rough_A \Smooth_A f^*(M)
    \cong \Rough_A f^* (\Smooth_B M)
    \\ \cong \Rough_A f^* (\Smooth_B \Rough_B M)
    \cong \Rough_A \Smooth_A f^* \Rough_B M
    \cong \Rough_A f^* (\Rough_B M).
  \end{multline*}
  Finally, adjoint associativity implies
  \begin{multline*}
    \Rough_A f^*(\Rough_B M)
    = \Hom_A(A,f^* \Hom_B(B,M))
    \\ \cong \Hom_B(A\hot_A B,M)
    \cong f^* \Hom_B(B,M) = f^*(\Rough_B M).
  \end{multline*}
  Putting things together, we obtain $\Rough_A f^*(M)\cong f^*(\Rough_B M)$ as
  desired.  Hence $f^!(M) = \Smooth_A f^*(\Rough_B M)\cong \Smooth_A
  f^*(M)\cong f^*(M)$ by Theorem~\ref{the:smoothening}.
\end{proof}

Evidently, $f_!$ and~$f^*$ commute with direct limits and $f_*$ and~$f^!$
commute with inverse limits.  Moreover, the functors $f^!$ and~$f^*$ are exact
because~$B$ is projective as a left or right $B$\nbd{}module.  Therefore,
their adjoint functors $f_!$ and~$f_*$ preserves projectives and injectives,
respectively.  This can also be seen directly from
$$
f_!(A\hot X)\cong B\hot X,
\qquad
f_*(\Smooth_A \Hom(A,X))\cong \Smooth_B \Hom(B,X).
$$
The first assertion is trivial, the second one follows from adjoint
associativity.

The same constructions define triangulated functors
$$
f_!,f_*\colon \Ho(A)\to\Ho(B),
\qquad
f^*,f^!\colon \Ho(B)\to\Ho(A).
$$
By exactness, $f^*,f^!\colon \Ho(B)\to\Ho(A)$ descend to functors
$f^*,f^!\colon \Der(B)\to\Der(A)$.  Of course, we still have $f^*=f^!$ on
$\Der(B)$ if~$f$ is proper.  The functors $f_!$ and~$f_*$ need not be exact.
As a substitute, we have their derived functors
$$
\Left f_!, \Right f_*\colon \Der(A)\to\Der(B).
$$
They are defined by applying~$f_!$ to a projective approximation and~$f_*$
to an injective approximation, respectively.  Thus
$$
\Left f_!(M) = B\Lhot_A M,
\qquad
\Right f_*(M) = \Smooth_B \bigl(\Right\Hom_A(B,M)\bigr).
$$

\begin{theorem}  \label{the:induced_adjoint_Der}
  For any morphism $f\colon A\dashrightarrow B$ there are natural isomorphisms
  \begin{align}  \label{eq:morphism_induced_ra}
    \Right\Hom_B(N,\Right f_*(M))\cong \Right\Hom_A(f^*(N),M),
  \\  \label{eq:morphism_induced_la}
    \Right \Hom_B(\Left f_!(M),N)\cong \Right\Hom_A(M, f^!(N)),
  \end{align}
  in~$\Ho$.  Thus $\Right f_*$ is right adjoint to $f^*\colon
  \Der(B)\to\Der(A)$ and $\Left f_!$ is left adjoint to $f^!\colon
  \Der(B)\to\Der(A)$.
\end{theorem}

\begin{proof}
  Observe first that the adjointness relations of
  Proposition~\ref{pro:adjointness_morphism} extend to the homotopy categories
  of chain complexes.  Let $P(M)\to M\to I(M)$ be the projective and injective
  approximation functors.  We have
  \begin{multline*}
    \Right\Hom_B(N,\Right f_*(M))
    \cong \Right\Hom_B(N,f_*(I(M)))
    \cong \Hom_B(N,f_*(I(M)))
    \\ \cong \Hom_A(f^*N,I(M))
    \cong \Right\Hom_A(f^*N,M)
  \end{multline*}
  because~$f_*$ preserves injectives.  This argument also uses
  Proposition~\ref{pro:preserve_projective}.  A similar computation, which
  uses that~$f_!$ preserves projectives,
  yields~\eqref{eq:morphism_induced_la}.
\end{proof}

There are canonical maps
\begin{align}
  \label{eq:fust_RightHom}
  f^*\colon \Right\Hom_B(M,N) \to \Right\Hom_A(f^*M,f^*N),
  \\
  \label{eq:fush_RightHom}
  f^!\colon \Right\Hom_B(M,N) \to \Right\Hom_A(f^!M,f^!N),
  \\
  \label{eq:flst_RightHom}
  f_*\colon \Right\Hom_A(M,N) \to \Right\Hom_B(\Right f_*M,\Right f_*N),
  \\
  \label{eq:flsh_RightHom}
  f_!\colon \Right\Hom_A(M,N) \to \Right\Hom_B(\Left f_!M, \Left f_!N)
\end{align}
in $\Ho$.  They have the characteristic property that when passing to
homology, we get the action of the functors $f^*$, $f^!$, $\Right f_*$, and
$\Left f_!$ on the spaces of morphisms in the derived categories.  There is a
natural map
\begin{equation}  \label{eq:induced_Lhot}
  f_*\colon f^*(M)\Lhot_A f^*(N) \to M\Lhot_B N
\end{equation}
as well.  The constructions are essentially the same in all cases.  Let us
construct~$f^*$.  Choose a projective approximation $P(M)\to M$ and consider
$f^*(P(M))\inOb\Ho(A)$.  Since~$f^*$ is exact, the map $f^*(P(M))\to f^*(M)$
is a quasi\brd{}isomorphism.  Hence a projective approximation for $f^*(P(M))$
is a projective approximation for $f^*(M)$ as well.  Hence we get bounded
linear chain maps
\begin{multline*}
  \Right\Hom_B(M,N)
  \cong \Hom_B(P(M),N)
  \to \Hom_A(f^* P(M),f^*N)
  \\ \to \Right\Hom_A(f^* P(M), f^*N)
  \cong \Right\Hom_A(f^* M, f^*N).
\end{multline*}
The resulting map $\Right\Hom_B(M,N)\to\Right\Hom_A(f^*M,f^*N)$ does not
depend on auxiliary choices and yields a natural transformation between these
bifunctors.  Moreover, the construction of~$f^*$ is functorial (up to chain
homotopy), that is, $\ID^*=\ID$ and $f^*g^*=(gf)^*$ in $\Ho$.

\section{Isocohomological morphisms}
\label{sec:isocohomological}

Let $A$ and~$B$ be two quasi-unital bornological algebras and let $i\colon
A\dashrightarrow B$ be a morphism.  We want to investigate when $i^*\colon
\Der(B)\to\Der(A)$ is fully faithful.

Consider first the functor $i^*\colon \Mod(B)\to\Mod(A)$.  Since $i^*M$
and~$M$ have the same underlying bornological vector space, we get a
bornological embedding $\Hom_A(i^*M,i^*N)\subseteq \Hom_B(M,N)$ for all
$M,N\inOb\Mod(B)$.  Equality holds, for instance, if~$i$ is a proper morphism
that has dense range in the sense that $\cl{i(A)}=B$.  Thus $i^*\colon
\Mod(B)\to\Mod(A)$ is always faithful, and fully faithful for proper morphisms
with dense range.  However, this does not suffice to conclude that $i^*\colon
\Der(B)\to\Der(A)$ is fully faithful because of the following proposition.

\begin{proposition}  \label{pro:isocoh_extensions}
  If $i^*\colon \Der(B)\to\Der(A)$ is fully faithful, then the subcategory
  $i^*\Mod(B)\subseteq \Mod(A)$ is closed under extensions with a bounded
  linear section.
\end{proposition}

\begin{proof}
  We can describe $\Ext^1_B(M,N)$ using either morphisms in the derived category or
  isomorphism classes of conflations.  The first description yields
  $$
  \Ext^1_B(M,N)\cong\Ext^1_A(i^*M,i^*N)
  \qquad \forall M,N\inOb\Mod(B)
  $$
  if $i^*\colon \Der(B)\to\Der(A)$ is fully faithful.  By the second
  description, this means that any conflation $M\into E\prto N$ in $\Mod(A)$
  with $M,N\inOb\Mod(B)$ is isomorphic to one in $i^*\Mod(B)$.  In
  particular, $E\in i^*\Mod(B)$, as desired.
\end{proof}

\begin{example}  \label{exa:integers_not_isocoh}
  Consider the Banach convolution algebras $A=\Sch^1(\Z)$ and $B=\ell_1(\Z)$
  on the group of integers~$\Z$, which are defined by the norms
  $$
  \norm{f}^1\defeq \sum_{n\in\Z} \abs{f(n)} (\abs{n}+1),
  \qquad
  \norm{f}^0\defeq \sum_{n\in\Z} \abs{f(n)}.
  $$
  Let $i\colon A\to B$ be the inclusion.  We claim that
  $i^*\Mod(B)\subseteq \Mod(A)$ is not closed under extensions.  Hence
  $i^*\colon \Der(B)\to\Der(A)$ cannot be fully faithful.  A similar
  counterexample exists for the embedding $C^1(S^1)\to C(S^1)$.

  Equip~$\C$ with the trivial representation of~$\Z$ and the corresponding
  module structure over $A$ and~$B$.  Consider the representation of~$\Z$
  on~$\C^2$ by
  $$
  \Z\ni n \mapsto \begin{pmatrix}  1 & n \\ 0 & 1  \end{pmatrix}.
  $$
  This is an extension of the trivial representation by itself.  Since the
  norm of these matrices grows linearly in~$n$, it defines a module over
  $\Sch^1(\Z)$, but not over $\ell_1(\Z)$.  Hence it is an element of
  $\Ext^1_A(\C,\C)$ that does not belong to $\Ext^1_B(\C,\C)$.
\end{example}

Recall that $\Ho=\Ho(\Born)$ is the homotopy category of chain complexes over
$\Born$.  The natural transformations $i^*\colon
\Right\Hom_B(M,N)\to\Right\Hom_A(i^*M,i^*N)$ and $i_*\colon i^*M\Lhot_A i^*N
\to M\Lhot_B N$ used below are constructed in \eqref{eq:fust_RightHom}
and~\eqref{eq:induced_Lhot}.

\begin{theorem}  \label{the:isocoh}
  Let $i\colon A\dashrightarrow B$ be a morphism between two quasi-unital
  algebras.  Let $P_\bullet\to A$ be a projective $A$\nbd{}bimodule resolution
  of $A\inOb\Mod(A\hot A^\op)$.  The following conditions are equivalent:
  \begin{enumerate}[\ref{the:isocoh}.1.]

  \item $B\hot_A P_\bullet \hot_A B$ is a projective $B$\nbd{}bimodule
    resolution of~$B$.

  \item The map $i_*\colon B\Lhot_A B\to B\Lhot_B B\cong B$ is an
    isomorphism in~$\Ho$.

  \item The map~$\alpha$ in the commuting diagram
    $$
    \xymatrix{
      (\Left i_!)\circ i^*(M) \ar@{=}[d] \ar[r]^-{\alpha} & M \\
      i^*(B)\Lhot_A i^*(M) \ar[r]^-{i_*} & B\Lhot_B M \ar[u]^{\cong}
    }
    $$
    is an isomorphism in $\Der(B)$ for all $M\inOb\Der(B)$.

  \item The map $i_*\colon i^*(M)\Lhot_A i^*(N) \to M\Lhot_B N$ is an
    isomorphism in~$\Ho$ for any $M\inOb\Der(B^\op)$, $N\inOb\Der(B)$.

  \item The map $i^*\colon \Right\Hom_B(B,\Hom(B,X)) \to \Right\Hom_A(B,
    \Hom(B,X))$ is an isomorphism in~$\Ho$ for all $X\inOb\Born$.

  \item The map $i^*\colon \Right\Hom_B(B,M) \to \Right\Hom_A(i^*B,i^*M)$
    is an isomorphism in~$\Ho$ for all $M\inOb\Der(B)$.

  \item The map $i^*\colon \Right\Hom_B(M,N) \to \Right\Hom_A(i^*M,i^*N)$
    is an isomorphism in~$\Ho$ for all $M,N\inOb\Der(B)$.

  \item The natural map
    $$
    M \cong \Smooth_B \Right\Hom_B(B,M)
    \to \Smooth_B \Right\Hom_A(B,M)
    = (\Right i_*)\circ i^*(M)
    $$
    is an isomorphism in $\Der(B)$ for all $M\inOb\Der(B)$.

  \item The functor $i^*\colon \Der(B)\to\Der(A)$ is fully faithful.

  \end{enumerate}
\end{theorem}

\begin{proof}
  First we check the equivalence of \ref{the:isocoh}.1 and \ref{the:isocoh}.2.
  The functor $X\mapsto X\hot_A B$ maps projective $A$\nbd{}bimodules to
  projective $A,B$\brd{}bimodules.  Similarly, $Y\mapsto B\hot_A Y$ maps
  projective $A,B$\brd{}bimodules to projective $B$\nbd{}bimodules.  Hence
  $P_\bullet\hot_A B$ is projective in $\Ho(A\hot B^\op)$ and $B\hot_A
  P_\bullet \hot_A B$ is projective in $\Ho(B\hot B^\op)$.  Since~$A$ is
  projective as a right $B$\nbd{}module, the map $P_\bullet\to A$ is an
  isomorphism in $\Ho(A^\op)$.  Thus $P_\bullet\hot_A B\to A\hot_A B\cong B$
  is a resolution.  It is a projective $A,B$\brd{}bimodule resolution of~$B$.
  It follows that $B\Lhot_A B\cong B\hot_A (P_\bullet\hot_A B)$.  This is a
  projective $B$\nbd{}bimodule resolution if and only if it is a resolution,
  if and only if the map $B\Lhot_A B\to B$ is an isomorphism in~$\Ho$.  Thus
  \ref{the:isocoh}.1 and \ref{the:isocoh}.2 are equivalent.
  
  We have $B\Lhot_B M\cong B\hot_B M \cong M$ because~$B$ is projective
  and~$M$ is essential.  Hence the map~$\alpha$ in the diagram in
  \ref{the:isocoh}.3 is an isomorphism if and only if~$i_*$ is.  Thus
  \ref{the:isocoh}.4 implies \ref{the:isocoh}.3 and \ref{the:isocoh}.3 implies
  \ref{the:isocoh}.2.  We claim that, conversely, \ref{the:isocoh}.2 implies
  \ref{the:isocoh}.4.  The class of $M\inOb\Der(B^\op)$ with $i^*(M)\Lhot_A
  i^*(N) \cong M\Lhot_B N$ for all $N\inOb\Der(B)$ is a triangulated
  subcategory of $\Der(B^\op)$ and also closed under direct sums.  By the
  quasi-unital variant of Proposition~\ref{pro:free_generate_DerA}, this
  subcategory is equal to $\Der(B^\op)$ once it contains $X\hot B$ for all
  $X\inOb\Born$.  Fix such an~$X$.  The class of $N\inOb\Der(B)$ with $(X\hot
  B)\Lhot_A i^*(N) \cong (X\hot B)\Lhot_B N$ is a triangulated subcategory of
  $\Der(B)$ closed under direct sums.  Again by
  Proposition~\ref{pro:free_generate_DerA}, it is equal to $\Der(B)$ once it
  contains $B\hot Y$ for all $Y\inOb\Born$.  Thus $i^*(M)\Lhot_A i^*(N) \cong
  M\Lhot_B N$ holds for all $M\inOb\Der(B^\op)$, $N\inOb\Der(B)$, once it
  holds for $M=X\hot B$ and $N=B\hot Y$.  Since $X\hot\blank$ is an exact
  functor for extensions with a bounded linear section, we have $(X\hot
  B)\Lhot_A (B\hot Y)\cong X\hot (B\Lhot_A B)\hot Y$ and $(X\hot B)\Lhot_B
  (B\hot Y)\cong X\hot (B\Lhot_B B)\hot Y \cong X\hot B\hot Y$.  Thus
  \ref{the:isocoh}.4 follows from \ref{the:isocoh}.2 as asserted.
  
  It is also clear that \ref{the:isocoh}.7 implies \ref{the:isocoh}.6 and
  \ref{the:isocoh}.6 implies \ref{the:isocoh}.5.  Almost the same argument
  that yields the implication
  \ref{the:isocoh}.2$\Longrightarrow$\ref{the:isocoh}.4 also shows that
  \ref{the:isocoh}.5 implies \ref{the:isocoh}.7.  Hence \ref{the:isocoh}.5--7
  are equivalent.
  
  Next we claim that \ref{the:isocoh}.2 and \ref{the:isocoh}.5--7 are
  equivalent.  Equation~\eqref{eq:derived_adjoint_associativity} yields
  \begin{align*}
    \Right\Hom_A(B,\Hom(B,X)) &\cong \Hom(B\Lhot_A B,X),
    \\
    \Right\Hom_B(B,\Hom(B,X)) &\cong \Hom(B\Lhot_B B,X),
  \end{align*}
  for all $X\inOb\Ho(\Born)$.  (We do not have to derive $\Hom(\blank,X)$
  because this functor is exact on extensions with a bounded linear section.)
  The map
  \begin{equation}
    \label{eq:silly}
    i^*\colon \Right\Hom_B(B,\Hom(B,X))\to\Right\Hom_A(B,\Hom(B,X))
  \end{equation}
  in \ref{the:isocoh}.5--7 corresponds to the map induced by $i_*\colon
  B\Lhot_A B\to B\Lhot_B B$ under these isomorphisms.  By the Yoneda Lemma,
  \eqref{eq:silly} is an isomorphism for all $X\inOb\Ho$ if and only if
  $i_*\colon B\Lhot_A B\to B\Lhot_B B$ is an isomorphism in~$\Ho$.  Hence
  \ref{the:isocoh}.2 implies \ref{the:isocoh}.5 and \ref{the:isocoh}.7 implies
  \ref{the:isocoh}.2.  So far, we have seen that \ref{the:isocoh}.1--7 are
  equivalent.
  
  We claim that \ref{the:isocoh}.6 and \ref{the:isocoh}.8 are equivalent.  It
  is clear that \ref{the:isocoh}.6 implies \ref{the:isocoh}.8 by applying the
  functor $\Smooth_B\colon \Der(B^+)\to\Der(B)$.  Conversely, we claim that we
  get \ref{the:isocoh}.6 back from \ref{the:isocoh}.8 if we apply
  $\Rough_B\colon \Der(B)\to\Der(B^+)$.  We have natural isomorphisms
  $$
  \Rough_B \Smooth_B \Right\Hom_B(B,M)
  \cong
  \Rough_B \Smooth_B \Rough_B M
  \cong
  \Rough_B M
  \cong \Right\Hom_B(B,M)
  $$
  because~$B$ is projective and
  $\Rough_B\circ\Smooth_B=\Rough_B^2=\Rough_B$.  Using adjoint associativity,
  we also get
  \begin{multline*}
    \Rough_B \Smooth_B \Right\Hom_A(B,M)
    \cong \Right\Hom_B(B,\Right\Hom_A(B,M))
    \\ \cong \Right\Hom_A(B\Lhot_B B,M)
    \cong \Right\Hom_A(B,M).
  \end{multline*}
  One checks easily that these isomorphisms are compatible with the maps~$i^*$
  on both sides.  Hence \ref{the:isocoh}.6$\iff$\ref{the:isocoh}.8.
  
  It is clear that \ref{the:isocoh}.7 implies \ref{the:isocoh}.9 by passing to
  homology.  For the converse, we claim that \ref{the:isocoh}.8 and
  \ref{the:isocoh}.9 are equivalent.  This follows in a formal way from the
  adjointness of $\Right i_*$ and~$i^*$.  It yields natural isomorphisms
  $\Der_A(i^*X,i^*M) \cong \Der_B(X,(\Right i_*)\circ i^*M)$ for all
  $X,M\inOb\Der(B)$.  The composite
  $$
  \Der_B(X,M)
  \overset{i^*}\to \Der_A(i^*X,i^*M)
  \congto \Der_B(X,(\Right i_*)\circ i^*M)
  $$
  is induced by some map $M\to (\Right i_*)\circ i^*M$, which turns out to
  be the one described in \ref{the:isocoh}.8.  The Yoneda Lemma implies that
  $\Der_B(X,M) \to \Der_A(i^*X,i^*M)$ is an isomorphism for all~$X$ if and
  only if $M\to(\Right i_*)\circ i^*M$ is an isomorphism.  This means that
  \ref{the:isocoh}.8 and \ref{the:isocoh}.9 are equivalent.
\end{proof}

\begin{definition}  \label{def:isocoh}
  A morphism between quasi-unital bornological algebras is called
  \emph{isocohomological} if the equivalent conditions of
  Theorem~\ref{the:isocoh} are satisfied.
\end{definition}

\begin{proposition}  \label{pro:compose_isocoh}
  Let $f\colon A\dashrightarrow B$ and $g\colon B\dashrightarrow C$ be
  composable morphisms between quasi-unital bornological algebras.  If two of
  $f$, $g$, $g\circ f$ are isocohomological, so is the third.
\end{proposition}

\begin{proof}
  Suppose first that~$g$ is isocohomological.  Then it follows immediately
  from condition~\ref{the:isocoh}.9 that $g\circ f$ is isocohomological if and
  only if~$f$ is isocohomological.  Now suppose that $g\circ f$ and~$f$ are
  isocohomological.  We want to show that~$g$ is isocohomological.  Since~$f$
  is isocohomological, we have $g^*(C) \Lhot_B g^*(C)\cong f^*g^*(C) \Lhot_A
  f^*g^*(C)$ by \ref{the:isocoh}.4.  Since $g\circ f$ is isocohomological, the
  latter is isomorphic to~$C$ by \ref{the:isocoh}.2.  This means that~$g$ is
  isocohomological.
\end{proof}

\begin{proposition}  \label{pro:tensor_isocoh}
  Let $A_1\dashrightarrow B_1$ and $A_2\dashrightarrow B_2$ be
  isocohomological morphisms.  Then the induced morphism $A_1\hot
  A_2\dashrightarrow B_1\hot B_2$ is isocohomological as well.
\end{proposition}

\begin{proof}
  Use $(B_1\hot B_2) \Lhot_{A_1\hot A_2} (B_1\hot B_2) \cong (B_1\Lhot_{A_1}
  B_2) \hot (B_2\Lhot_{A_2} B_2)$.
\end{proof}

\section{Isocohomological group convolution algebras}
\label{sec:isocoh_groups}

Let $\CCINF(G)$ be the convolution algebra of smooth functions of compact
support on a locally compact group~$G$.  This is a quasi-unital bornological
algebra, and $\Mod(\CCINF(G))$ is isomorphic to the category $\Mod(G)$ of
smooth representations of~$G$ (see~\cite{Meyer:Smooth}).  We frequently
replace $\CCINF(G)$ by~$G$ in our notation, writing $\Lhot_G$ and~$\hot_G$ for
$\Lhot_{\CCINF(G)}$ and $\hot_{\CCINF(G)}$ and $\Right\Hom_G$ and $\Hom_G$ for
$\Right\Hom_{\CCINF(G)}$ and $\Hom_{\CCINF(G)}$.

A \emph{smooth convolution algebra on~$G$} is a bornological algebra $\T(G)$
of functions on~$G$ that contains $\CCINF(G)$ as a dense subalgebra, such that
the left and right regular representations of~$G$ on $\T(G)$ are smooth.
Equivalently, $\T(G)$ is an essential bimodule over $\CCINF(G)$.  Let $i\colon
\CCINF(G)\to\T(G)$ be the embedding.  It follows from
Lemma~\ref{lem:embed_quasiunital} that $\T(G)$ is always quasi-unital and
that~$i$ is a proper morphism.
Proposition~\ref{pro:adjointness_proper_morphism} asserts, among other things,
that
$$
i^* \Smooth_{\T(G)}=i^* \Smooth_{\CCINF(G)},
\qquad
i^* \Rough_{\T(G)}=i^* \Rough_{\CCINF(G)}.
$$
A $\T(G)$\brd{}module~$M$ is essential if and only if $i^*M$ is essential
if and only if the module structure is the integrated form of a smooth
representation of~$G$ (\cite{Meyer:Smooth}).

Since~$i$ has dense range, the functor $i^*\colon \Mod(\T(G))\to\Mod(G)$ is
fully faithful, that is, $\Hom_{\T(G)}(M,N) = \Hom_G(M,N)$ for any
$M,N\inOb\Mod(\T(G))$.  Thus $\Mod(\T(G))$ is a full subcategory of the
category $\Mod(G)$ of smooth representations of~$G$.  We call such
representations \emph{$\T(G)$\brd{}tempered}.  A smooth representation is
$\T(G)$\brd{}tempered if and only if its integrated form $\CCINF(G)\to\End(M)$
extends to a bounded homomorphism $\T(G)\to\End(M)$.

The Schwartz algebras of Abelian locally compact groups and the Schwartz
algebras of reductive groups over local fields defined by Harish-Chandra are
examples of smooth convolution algebras.  If~$G$ is discrete, then the
smoothness condition is empty and a smooth convolution algebra is just an
algebra that contains $\C[G]$ as a dense subalgebra.  For instance, the
algebra $\Sch(G)$ defined in~\eqref{eq:def_Schwartz_space} is a smooth
convolution algebra on~$G$.

\begin{definition}[\cite{Meyer:Primes_Rep}]  \label{def:iscoh_convolution}
  We call a smooth convolution algebra $\T(G)$ \emph{isocohomological} if the
  embedding $\CCINF(G)\to\T(G)$ isocohomological.
\end{definition}

\begin{definition}  \label{def:symmetric_convolution}
  A smooth convolution algebra $\T(G)$ on a locally compact group~$G$ is
  called \emph{symmetric} if $U\phi(g,h)\defeq \phi(gh,g)$ and its inverse
  $U^{-1}\phi(g,h)\defeq \phi(h,h^{-1}g)$ are bounded operators on
  $\T(G)\hot\T(G)$.
\end{definition}

Of course, $\CCINF(G)$ is always symmetric.

Let $\T(G)$ be a symmetric convolution algebra and let $\mu\colon
\T(G)\hot\T(G)\to\T(G)$ be the convolution.  Then $(\mu\circ U)\phi(g) =
\int_G \phi(g,h) \,dh$.  Hence the trivial representation $\phi\mapsto \int_G
\phi(h)\,dh$ is a module over $\T(G)$.

The convolution algebras $\Hol(G)$ and $\Sch^k(G)$ for
$k\in\R_+\cup\{\infty,\omega\}$ defined in~\cite{Meyer:Combable_poly} for a
finitely generated discrete group are easily seen to be symmetric.  The
Schwartz algebra of a locally compact Abelian group is symmetric as well
because the map~$U$ is associated to a group homomorphism in this case.  In
contrast, the Harish-Chandra-Schwartz algebra of a reductive group over a
local field cannot be symmetric because the trivial representation of such a
group is not tempered.  For the same reason, the Jolissaint algebra of a
discrete group with rapid decay is only symmetric for groups of polynomial
growth.

\begin{proposition}  \label{pro:symmetric_isocoh}
  A symmetric convolution algebra $\T(G)$ is isocohomological if and only if
  $\T(G)\Lhot_G \C\cong \T(G)\Lhot_{\T(G)} \C\cong \C$.  By definition,
  $\T(G)\Lhot_G \C \cong \T(G) \hot_G P_\bullet$ for any projective resolution
  $P_\bullet\to\C$ of the trivial representation of~$G$.
\end{proposition}

\begin{proof}
  For $M\inOb\Mod(G^\op)$, $N\inOb\Mod(G)$, equip $M\hot N$ with the inner
  conjugation action $(m\otimes n)\cdot g\defeq (m\cdot g\otimes g^{-1}\cdot
  n)$ for all $g\in G$, $m\in M$, $n\in N$, and the associated right
  $\CCINF(G)$\brd{}module structure.  Then we have a natural isomorphism
  $M\hot_G N\cong (M\hot N) \hot_G \C$ (see~\cite{Meyer:Smooth}).  This easily
  implies the corresponding statement for chain complexes $M$ and~$N$.  It is
  checked in~\cite{Meyer:Smooth} that $M\hot N$ is a projective
  $\CCINF(G)$\brd{}module if $M$ or~$N$ is projective.  Hence
  $$
  M \Lhot_G N \cong (M\hot N) \Lhot_G \C
  $$
  for all $M\inOb\Der(G^\op)$, $N\inOb\Der(G)$.  In particular,
  $$
  \T(G)\Lhot_G \T(G) \cong
  \bigl(\T(G)\hot\T(G)\bigr)^{\rho\hot\lambda}\Lhot_G \C,
  $$
  where the superscript $\rho\hot\lambda$ indicates that we equip
  $\T(G)\hot\T(G)$ with the inner conjugation representation of~$G$.
  
  The operator~$U$ of Definition~\ref{def:symmetric_convolution} intertwines
  the inner conjugation action $\rho\hot\lambda$ and the regular
  representation on the second tensor factor $1\hot\rho$.  Hence
  $$
  \bigl(\T(G)\hot\T(G)\bigr)^{\rho\hot\lambda} \cong
  \bigl(\T(G)\hot\T(G)\bigr)^{1\hot\rho}
  $$
  for symmetric convolution algebras.  Therefore,
  $$
  \T(G)\Lhot_G \T(G) \cong \bigl(\T(G)\hot\T(G)\bigr)^{1\hot\rho} \Lhot_G
  \C \cong \T(G)\hot(\T(G)\Lhot_G \C).
  $$
  This implies the assertion.
\end{proof}

\begin{theorem}  \label{the:isocoh_coarse}
  Let~$G$ be a finitely generated discrete group.  If~$G$ admits a combing of
  polynomial growth, then $\Sch(G)$ is isocohomological; if~$G$ admits a
  combing of subexponential growth, then $\Sch^\omega(G)$ is isocohomological;
  and if~$G$ admits a combing of exponential growth, then $\Hol(G)$ is
  isocohomological.
  
  In general, if $\T(G)$ is $\Hol(G)$, $\Sch^\omega(G)$, or $\Sch(G)$, then
  $\T(G)$ is isocohomological if and only if the chain complex that is denoted
  by $\T\tilde{C}_\bullet(G)$ in~\cite{Meyer:Combable_poly} has a bounded
  contracting homotopy.  This property is invariant under quasi-isometry.
\end{theorem}

\begin{proof}
  Let~$G$ be a finitely generated discrete group and let $\T(G)$ be one of
  the convolution algebras $\Hol(G)$ and $\Sch^k(G)$ for
  $k\in\R_+\cup\{\infty,\omega\}$.  It is easy to see that $\T(G)$ is
  symmetric.  For an appropriate choice of resolution $P_\bullet$, the complex
  $\T(G) \hot_G P_\bullet$ is exactly the complex that is denoted by $\T
  C_\bullet(G)$ in~\cite{Meyer:Combable_poly}.
  Proposition~\ref{pro:symmetric_isocoh} therefore yields that $\T(G)$ is
  isocohomological if and only if the chain complex
  $\T\tilde{C}_\bullet(G)$ is contractible.  It is shown
  in~\cite{Meyer:Combable_poly} that this condition is invariant under
  quasi-isometry and is satisfied in the presence of suitable combings.
\end{proof}

In particular, the above theorem applies to finitely generated free
non-Abelian groups, to hyperbolic groups, to automatic groups, and to finitely
generated Abelian groups.  All these classes of groups have combings of linear
growth.  The following is proven in~\cite{Meyer:Primes_Rep}:

\begin{theorem}  \label{the:isocoh_Sch_Abelian}
  The Bruhat-Schwartz algebra $\Sch(G)$ is isocohomological for any Abelian
  locally compact group~$G$.
\end{theorem}

A \emph{weight function} on~$G$ is a function $w\colon G\to\R_{>0}$ such that
$w(gh)\le w(g)w(h)$ for all $g,h\in G$.  We do not require any relationship
between $w(g)$ and $w(g^{-1})$.  If~$w$ is a weight function, then
$$
\ell_1(G,w) \defeq \biggl\{f\colon G\to\C \biggm\vert
\sum_{g\in G} \abs{f(g)}w(g)<\infty \biggr\}
$$
is a Banach algebra with respect to convolution.  The functions $w(g)\cdot
(\ell(g)+1)^k$ for $k\in\N$ and $w(g)\cdot \alpha^{\ell(g)}$ for $\alpha>1$
are also weight functions if~$w$ is.  We let
$$
\Sch(G,w)\defeq \bigcap_{k\in\N} \ell_1(G,w\cdot (\ell+1)^k).
$$
A subset of $\Sch(G,w)$ is bounded if and only if it is bounded in
$\ell_1(G,w\cdot (\ell+1)^k)$ for all $k\in\N$.  Thus $\Sch(G,w)$ is a
Fréchet-Schwartz space.  More generally, if~$W$ is a set of weight functions,
we let
$$
\Sch(G,W)\defeq \bigcap_{w\in W} \Sch(G,w).
$$
A subset of $\Sch(G,W)$ is bounded if it is bounded in $\Sch(G,w)$ for all
$w\in W$.  If~$W$ is countable, then $\Sch(G,w)$ is a Fréchet-Schwartz space.
For instance, if $W=\{n^\ell\mid n\in\N_{\ge2}\}$, then $\Sch(G,W)=\Hol(G)$
because the exponential growth of $(n+1)^\ell/n^\ell$ dominates the polynomial
growth of $(\ell+1)^k$.  The convolution algebras $\Sch(G,W)$ are usually not
symmetric.  Nevertheless, we can treat them as in~\cite{Meyer:Combable_poly}.

\begin{definition}  \label{def:weight_combing}
  Let~$G$ be a discrete group and let $f_j\colon G\to G$, $j\in\N$, be a
  combing of~$G$ as defined in~\cite{Meyer:Combable_poly}.  A set of
  weights~$W$ is called \emph{compatible} with the combing if it has the
  following property: for any $w\in W$ there is a finite linear combination
  $\bar{w} \defeq \sum_{w\in W} \alpha_w w$ such that
  $$
  w(f_j(g))\cdot w(f_j(g)^{-1}h) \le \bar{w}(g)\cdot \bar{w}(g^{-1}h)
  $$
  for all $g,h\in G$, $n\in\N$.
\end{definition}

\begin{theorem}  \label{the:isocoh_weighted}
  Let~$G$ be a discrete group and let $(f_j)_{j\in\N}$ be a combing of~$G$ of
  polynomial growth.  Let~$W$ be a set of weight functions on~$G$ that is
  compatible with the combing.  Then $\Sch(G,W)$ is isocohomological.
\end{theorem}

\begin{proof}
  The proof is very similar to the proof of the corresponding statement
  without weights in~\cite{Meyer:Combable_poly}.  Hence we only outline the
  necessary changes and assume that the reader is familiar
  with~\cite{Meyer:Combable_poly}.

  Let $V\defeq \Sch(G,W)\hot\Sch(G,W)$ equipped with the inner conjugation
  action.  We have to show that the natural chain map $V\Lhot_G \C \to
  \Sch(G,W)$ is an isomorphism in $\Ho$.  The operator~$U$ in
  Definition~\ref{def:symmetric_convolution} is a bornological isomorphism
  between~$V$ and
  \begin{multline*}
    \tilde{V} \defeq \{f\colon G×G\to\C \mid
    \\ \sum_{g,h\in G} \abs{f(g,h)} w(h) w(h^{-1}g) (\ell(g)+\ell(h)+1)^k
    <\infty \ \forall w\in W, k\in\N\}.
  \end{multline*}
  The isomorphism is equivariant for the regular representation $f\cdot
  g(x,y)\defeq f(x,y\cdot g)$ on~$\tilde{V}$.  The free chain complex
  $C_\bullet(G)$ can be viewed as a free $\C[G]$\brd{}module resolution
  of~$\C$.  Hence $V\Lhot_G \C \cong V\hot_G C_\bullet(G)$.  We identify
  $V\hot_G C_n(G)$ with the space of functions $f\colon G×G^{n+1}\to \C$ that
  satisfy the control condition in the definition of $\Sch C_\bullet(G)$
  in~\cite{Meyer:Combable_poly} with respect to the last $n+1$ variables and
  for which the function $(g,h)\mapsto f(g,h\cdot x)$ belongs to~$\tilde{V}$
  for all $x\in G^{n+1}$.  A subset of $V\hot_G C_n(G)$ is bounded if it
  satisfies the control condition uniformly and if the set of functions
  $(g,h)\mapsto f(g,h\cdot x)$ is bounded in~$\tilde{V}$ for any fixed $x\in
  G^{n+1}$.  The boundary map on $V\hot_G C_\bullet(G)$ is $\ID\otimes\delta$.
  
  We have $V\hot_G C_0(G)=\tilde{V}$.  The multiplication map $V\to\Sch(G,W)$
  corresponds to the map $\alpha\colon \tilde{V}\to \Sch(G,W)$, $\alpha
  f(g)=\sum_h f(g,h)$.  The map~$\alpha$ is split surjective.  Let~$K_\bullet$
  be the complex with $K_n=V\hot_G C_n(G)$ for $n\ge1$ and $K_0=\ker\alpha$.
  Instead of showing that~$\alpha$ is a homotopy equivalence, we may show
  that~$K_\bullet$ is contractible.  Let $H\colon
  \tilde{C}_\bullet(G)\to\tilde{C}_\bullet(G)$ be the contracting homotopy
  that is used in the proof of \cite{Meyer:Combable_poly}*{Theorem 4}.  Then
  $\ID\otimes H$ is a contracting homotopy for the dense subcomplex of
  compactly supported functions in~$K_\bullet$.  As
  in~\cite{Meyer:Combable_poly}, one checks that this extends to a bounded
  contracting homotopy for~$K_\bullet$.  The compatibility between the set of
  weights and the combing is exactly what is needed to prove that $\ID\otimes
  H$ is bounded with respect to the growth condition that defines~$\tilde{V}$.
\end{proof}

\begin{example}  \label{exa:group_hom_compatible}
  If $w\colon G\to\R_{>0}$ is a group homomorphism, then
  $$
  w(f_j(g))\cdot w(f_j(g)^{-1}h) = w(h) = w(g)\cdot w(g^{-1}h),
  $$
  no matter what combing we use.  Hence the compatibility is automatic
  if~$W$ is a set of group homomorphisms.  For instance, let $G=\Z^n$.  This
  group has a well-known combing of linear growth.  It is defined by
  approximating the straight line in~$\R^n$ by a path in~$\Z^n$.  Group
  homomorphisms $G\to\R_{>0}$ are of the form $x\mapsto \exp(a\cdot x)$ for
  some $a\in\R^n$.  Hence any subset $W\subseteq\R^n$ defines a weighted
  Schwartz algebra $\Sch(\Z^n,W)$.  All these convolution algebras are
  isocohomological by Theorem~\ref{the:isocoh_weighted}.
  
  By the way, $\max\{\exp(a\cdot x),\exp(b\cdot x)\}$ dominates $\exp(c\cdot
  x)$ if~$c$ is a convex combination of $a$ and~$b$.  Therefore,
  $\Sch(\Z^n,W)=\Sch(\Z^n,\bar{W})$ if~$\bar{W}$ is the convex hull of~$W$.
  Thus it suffices to consider the algebras $\Sch(\Z^n,W)$ for convex
  $W\subseteq\R^n$.  Any such convex subset is a union of an increasing
  sequence of convex polyhedra.  This implies that $\Sch(\Z^n,W)$ is always a
  Fréchet-Schwartz algebra.
\end{example}

\begin{example}  \label{exa:Jolissaint}
  Let~$\F_r$ be the free group on~$r$ generators $s_1,\dotsc,s_r$.  This group
  has the rapid decay property of Paul Jolissaint (\cite{Jolissaint:RD}).
  This means that
  $$
  \Sch_2(\F_r) \defeq \bigcap_{k\in\N} \ell^2(\F_r,(\ell+1)^k)
  $$
  is an algebra with respect to convolution.  Using the standard free
  $\C[\F_r]$\brd{}resolution of the trivial representation of length~$1$, we
  identify $(\Sch_2(\F_r)\hot\Sch_2(\F_r))^{\rho\hot\lambda}\Lhot_{\F_r} \C$
  with the complex $(\Sch_2(\F_r)\hot\Sch_2(\F_r))^r \overset{\delta}\to
  \Sch_2(\F_r)\hot\Sch_2(\F_r)$ of length~$1$, where
  $$
  \delta((x_j\otimes y_j))
  \defeq \sum_{j=1}^r x_j\otimes y_j - x_js_j^{-1}\otimes s_jy_j.
  $$
  One can show that~$\delta$ is not an isomorphism onto the kernel of the
  multiplication map $\Sch_2(\F_r)\hot\Sch_2(\F_r)\to\Sch_2(\F_r)$.  This
  means that $\Sch_2(\F_r)$ is \emph{not} isocohomological.
\end{example}

\subsection{Variants of the de Rham complex}
\label{sec:deRham_complexes}

The complexes that we introduce in this section are used to prove that certain
convolution algebras are isocohomological.

Let~$M$ be a smooth oriented manifold of dimension~$n$.  We will only use the
case $M=\R^n$ later.  However, since we equip~$\R^n$ with group actions and
want to do constructions equivariantly, it is useful to formulate some
definitions more generally.  Let $T^*M\to M$ be the cotangent bundle of~$M$
and let $\Lambda^*(M)\to M$ be its exterior algebra bundle.  This is a graded
algebra.  Let $\Omega(M)$ be the space of smooth sections of $\Lambda^*M$,
equipped with its usual Fréchet topology and the precompact bornology that it
defines.  The de Rham boundary~$d^\dR$ turns $\Omega(M)$ into a bornological
chain complex, whose homology is known as the de Rham cohomology of~$M$.

Let $\CCINF\Omega(M) \subseteq \Omega(M)$ be the subspace of compactly
supported smooth sections.  A subset of $\CCINF\Omega(M)$ is bounded if it is
bounded in $\Omega(M)$ and all all its elements have uniformly compact
support.  The orientation of~$M$ provides a bounded linear functional
$\smallint\colon \CCINF\Omega^n(M)\to\C$ that satisfies $\smallint\circ
d^\dR=0$.  Hence we get a chain complex
\begin{equation}  \label{eq:deRham_complex}
  0 \longrightarrow \CCINF\Omega^0(M)
  \overset{d^\dR}\longrightarrow \CCINF\Omega^1(M)
  \overset{d^\dR}\longrightarrow \dotso
  \overset{d^\dR}\longrightarrow \CCINF\Omega^n(M)
  \overset{\smallint}\longrightarrow \C.
\end{equation}
This chain complex is natural for orientation preserving diffeomorphisms.  If
a diffeomorphism $\Phi\colon M\to M$ reverses the orientation, we let it act
on $\CCINF\Omega^j(M)$ by $\omega\mapsto -\Phi^*\omega$ for all~$j$.  Thus
$\smallint$ and~$d^\dR$ commute with the action of~$\Phi$.
Thus~\eqref{eq:deRham_complex} becomes a chain complex over $\Mod(G)$ if~$G$
is a Lie group that acts smoothly on~$M$.  We equip~$\C$ with the trivial
representation of~$G$.  Since we want~\eqref{eq:deRham_complex} to be a
resolution, we shift the grading so that $\CCINF\Omega^n(M)$ occurs in degree
$n-j$.

\begin{proposition}  \label{pro:deRham_resolution}
  Suppose that a Lie group~$G$ (possibly with infinitely many components) acts
  properly and smoothly on~$\R^n$ for some $n\in\N$.  Then the chain
  complex~\eqref{eq:deRham_complex} for $M=\R^n$ is a projective resolution of
  the trivial representation of~$G$.
\end{proposition}

\begin{proof}
  We have to check that $\CCINF\Omega(\R^n)$ is a projective
  $\CCINF(G)$\nbd{}module and that the complex~\eqref{eq:deRham_complex} for
  $M=\R^n$ has a bounded contracting homotopy.  The second assertion is
  well-known.  Using the isomorphism $\CCINF\Omega(\R^n) \cong
  \CCINF\Omega(\R)^{\hot n}$, we can reduce the proof to the easy case $n=1$.
  To prove projectivity of $\CCINF\Omega(\R^n)$, we need a
  $G$\nbd{}equivariant section for the map
  $$
  \alpha\colon \CCINF(G\times M,\Lambda^*M) \congto
  \CCINF(G)\hot\CCINF\Omega(M) \to
  \CCINF\Omega(M)
  $$
  that defines the group action.  There is $\psi\in\CCINF(M)$ with $\int_G
  \psi(g^{-1}x)\,d\lambda(g)=1$ for all $x\in M$.  The map $\sigma
  f(g,x)\defeq \psi(x)\cdot g^{-1}_*f(gx)$ is the desired section
  for~$\alpha$.
\end{proof}

Next we describe variants $\Sch\Omega(\R^n)$, $\Sch^\omega\Omega(\R^n)$, and
$\Hol(\R^n)$ of~\eqref{eq:deRham_complex} for $M=\R^n$.  We let
$\Sch\Omega(\R^n)$ be the space of all sections of $\Lambda^*\R^n$ with rapid
decay; that is, for all $k\in\N$ and all constant coefficient differential
operators~$D$ on~$\R^n$, there is $C>0$ with $\norm{D(f)(x)}\le C\cdot
(\norm{x}+1)^{-k}$.  Here we use the usual Euclidean metric on~$\R^n$ to
define the norm on the fibres of $\Lambda^*\R^n$.  A subset of
$\Sch\Omega(\R^n)$ is bounded if the above estimates hold uniformly for its
elements.  The operators $d^\dR$ and $\smallint$ extend to $\Sch\Omega(\R^n)$,
so that $\Sch\Omega(\R^n)$ becomes a chain complex of bornological vector
spaces.

Similarly, we define $\Hol\Omega(\R^n)$ by requiring exponential decay, that
is, for all $\alpha>1$ and all constant coefficient differential operators~$D$
on~$\R^n$, there is $C>0$ such that $\norm{D(f)(x)}\le C\cdot
\alpha^{\norm{x}}$.  Finally, we define $\Sch^\omega\Omega(\R^n)$ by a
subexponential decay condition: there is $\alpha>1$ such that for all constant
coefficient differential operators~$D$ on~$\R^n$ there is $C>0$ with
$\norm{D(f)(x)}\le C\cdot \alpha^{\norm{x}}$.  In both cases, we call a subset
bounded if its elements satisfy the appropriate estimate uniformly.

\begin{lemma}  \label{lem:tempered_deRham_complex}
  The complexes $\Sch\Omega(\R^n)$, $\Sch^\omega\Omega(\R^n)$ and
  $\Hol\Omega(\R^n)$ (all augmented by~$\smallint$) are contractible.
\end{lemma}

\begin{proof}
  One checks easily that $\Sch\Omega(\R^n)\cong \Sch\Omega(\R)^{\hot n}$.
  Hence the assertion for general~$n$ reduces to the special case $n=1$, which
  is easy.  The same argument works for $\Sch^\omega\Omega(\R^n)$ and
  $\Hol\Omega(\R^n)$.
\end{proof}

\subsection{Application to groups of polynomial growth}
\label{sec:pol_growth}

Let~$G$ be a discrete group of polynomial growth.  That is, $G$ is finitely
generated and the number of elements $g\in G$ satisfying $\ell(g)\le R$ has
polynomial growth as a function of~$R$.  This section is devoted to the proof
of the following theorem:

\begin{theorem}  \label{the:isocoh_poly_growth}
  Let~$G$ be a discrete group of polynomial growth.  Then the convolution
  algebras $\Sch(G)$, $\Sch^\omega(G)$ and $\Hol(G)$ on~$G$ are
  isocohomological.
\end{theorem}

Let $\T(G)$ be one of the convolution algebras in the statement of the
theorem.  We shall use the following known structure theorem:

\begin{theorem}  \label{the:poly_growth}
  Let~$G$ be a finitely generated discrete group of polynomial growth.  Then
  there exists a connected nilpotent Lie group~$\bar{G}$ and a cocompact
  lattice $G'\subseteq \bar{G}$ that is quasi-isometric to~$G$.
\end{theorem}

\begin{proof}
  By a celebrated result of Mikhail Gromov (\cite{Gromov:PolyNilpotent}), the
  group~$G$ contains a nilpotent subgroup~$G_1$ of finite index.  It follows
  from results of~\cite{Hall:Nilpotent} that~$G_1$ and hence~$G$ contain a
  torsion free subgroup~$G'$ of finite index, which is, of course, again
  finitely generated and nilpotent.  A famous result of Anatoli\u\i{}
  Ivanovi\v{c} Malcev (\cites{Malcev:Nil_homogeneous,
    Malcev:Nil_homogeneous_translation}) yields that~$G'$ is isomorphic to a
  cocompact lattice in a connected nilpotent Lie group.
\end{proof}

It follows from~\cite{Meyer:Combable_poly} that the assertions in
Theorem~\ref{the:isocoh_poly_growth} are invariant under quasi\brd{}isometry.
Hence we may assume in the following that~$G$ is a cocompact lattice in a
connected nilpotent Lie group~$\bar{G}$.  Let~$\LG$ be the Lie algebra
of~$\bar{G}$.  Since~$\bar{G}$ is nilpotent, the exponential map $\exp\colon
\LG\to\bar{G}$ is a diffeomorphism, that is, its inverse $\log\colon
\bar{G}\to\LG$ is everywhere defined and smooth.  In the following, we
identify $\bar{G}\cong\LG\cong\R^n$ and equip~$\LG$ with the action of~$G$ by
left translation.  This is a smooth action $G\times\LG\to\LG$ as in
Proposition~\ref{pro:deRham_resolution}.  Hence the complex
$\CCINF\Omega(\bar{G})$ is a projective resolution of the trivial
representation of~$G$.

\begin{lemma}  \label{lem:basic_nilpotent}
  $\T(G)\hot_G \CCINF\Omega(\bar{G}) \cong \T\Omega(\LG)$ as bornological
  chain complexes.
\end{lemma}

Lemma~\ref{lem:basic_nilpotent} and Lemma~\ref{lem:tempered_deRham_complex}
yield $\T(G)\Lhot_G \C\cong\C$.  Since $\T(G)$ is a symmetric convolution
algebra, this implies that $\T(G)$ is isocohomological and finishes the proof
of Theorem~\ref{the:poly_growth}.

\begin{proof}
  Let $X_1,\dotsc,X_n$ be an orthonormal basis for~$\LG$.  We view
  $X_1,\dots,X_n$ as left invariant vector fields on~$\bar{G}$ and then pull
  them back to~$\LG$ using the exponential map.  We still write
  $X_1,\dotsc,X_n$ for the resulting vector fields on~$\LG$.  Alternatively,
  we may also extend $X_1,\dots,X_n$ to translation invariant vector fields on
  the Euclidean space~$\LG$.  We denote these by $Y_1,\dotsc,Y_n$.  Since they
  are both $\CCINF(\LG)$\brd{}module bases for the space of vector fields
  on~$\LG$, they are related by pointwise multiplication with some function
  $\LG\to\Gl(\LG)$.  This function can be described as follows.  Let
  $\ad\colon\LG\to\End(\LG)$ be the adjoint representation of~$\LG$.  We need
  the holomorphic function
  $$
  \sigma\colon \C\to\C,\qquad \sigma(x)\defeq \frac{\sinh(x/2)}{x/2}.
  $$
  It is shown in \cite{Mostow:Strong_Rigidity} and~\cite{Bridson-Haefliger}
  that $Y_j(Z)$ is obtained from $X_j(Z)$ by applying the linear operator
  $\sigma\bigl(\ad(Z)\bigr)$ for all $Z\in\LG$.  Since~$\bar{G}$ is nilpotent,
  the function $\sigma(\ad Z)$ and its inverse $\sigma(\ad Z)^{-1}$ are both
  polynomial functions on~$\LG$.  Therefore, if $\phi\in\CINF(\LG)$, then
  $P(X_1,\dotsc,X_n)(\phi)$ has rapid decay for all polynomials~$P$ if and
  only if $P(Y_1,\dotsc,Y_n)(\phi)$ has rapid decay for all polynomials~$P$.
  
  Equip~$\bar{G}$ with some Riemannian metric that is invariant under left
  multiplication and let $d\colon \bar{G}\times\bar{G}\to\R_+$ be the
  associated metric.  It is well-known that the curves $\exp(tX)$ for
  $X\in\LG$ with $\norm{X}=1$ are unit speed geodesics in~$\bar{G}$.  Thus
  $d(x,1)=\norm{\log(x)}$ for all $x\in\bar{G}$.  Since $G\subseteq\bar{G}$ is
  cocompact, the word length function on~$G$ is equivalent to the length
  function $\norm{\log(x)}$, so that both norms give the same rapid decay
  condition.  This allows to identify $\Sch(G)\hot_G\CCINF(\LG)$ with the
  space of functions $\LG\to\C$ for which $P(X_1,\dotsc,X_n)(\phi)$ has rapid
  decay for all polynomials~$P$.  Since we may replace $(X_j)$ by $(Y_j)$,
  this gives $\Sch(\LG)$ as desired.  The same argument works for
  $\Sch^\omega(\bar{G})$ and $\Hol(\bar{G})$.
  
  Let $dX_1,\dots,dX_n$ be the basis for the space of covector fields that is
  dual to $X_1,\dots,X_n$.  Then the standard basis vectors
  $dX_{i_1}\wedge\dotsb\wedge dX_{i_j}$ for $\Lambda^*\LG$ provide a
  $\bar{G}$\nbd{}invariant basis for the space of differential forms on~$\LG$.
  This yields a $\bar{G}$\nbd{}equivariant isomorphism $\CCINF(\bar{G})\hot
  \Lambda^*\LG \cong \CCINF\Omega(\bar{G})$ and hence $\T(G)\hot_G
  \CCINF\Omega(\bar{G}) \cong \T(\LG)\hot\Lambda^*\LG$.  Here we use the
  differential forms $dX_{i_1}\wedge\dotsb\wedge dX_{i_j}$ as our basis to
  measure the decay of a differential form.  Since we get the same decay
  condition if we use $dY_{i_1}\wedge\dotsb\wedge dY_{i_j}$ instead, this
  space is equal to $\T\Omega(\LG)$ as desired.
\end{proof}

\subsection{Some smooth convolution algebras on Lie groups}
\label{sec:nilpotent_reductive_Lie}

The techniques used to prove Theorem~\ref{the:poly_growth} also apply to
certain smooth convolution algebras on Lie groups.  First we define the
convolution algebras that we are interested in.  Let~$G$ be a Lie group with
finitely many connected components.  (More generally, similar constructions
work for arbitrary almost connected locally compact groups).  Let~$d\lambda$
be a left invariant Haar measure on~$G$.  We need the following analogue of a
word length function.  There is a compact subset $S\subseteq G$ with
$S=S^{-1}$ and $\bigcup S^n=G$.  Define $\ell(g)$ to be the minimal $n\in\N$
with $g\in S^n$.  It is easy to see that the norms
$$
\norm{\phi}^k \defeq \int_G \abs{\phi(g)} (\ell(g)+1)^k\,d\lambda(g),
\qquad
\enorm{\phi}_\alpha \defeq \int_G \abs{\phi(g)} \alpha^{\ell(g)}\,d\lambda(g)
$$
are submultiplicative and hence define Banach convolution algebras
$\ell_1(G,(\ell+1)^k)$ and $\ell_1(G,\alpha^\ell)$ on~$G$.

Consider the space $A(G)\defeq \bigcap_{\alpha>1} \ell_1(G,\alpha^\ell)$ of
functions of exponential decay, equipped with the obvious bornology.  This is
the von Neumann bornology on a Fréchet space.  We let $\Hol(G)$ be the
smoothening of the left regular representation of~$G$ on $A(G)$.  Let
$\UE(\LG)$ be the universal enveloping algebra of~$G$, identified with the
space of distributions on~$G$ supported at the identity element.  Since $A(G)$
is bornologically metrisable, the smoothening $\Hol(G)$ is the subspace of
$\phi\in A(G)$ for which $D\ast\phi\in A(G)$ for all $D\in\UE(\LG)$.  A
subset~$S$ of $A(G)$ is bounded if and only if $D\ast S$ is bounded in $A(G)$
for each $D\in\UE(\LG)$ (see~\cite{Meyer:Smooth}).  If~$G$ is only almost
connected, the smoothening is a bit more complicated to describe explicitly.
We assume in the following that~$G$ is a Lie group to simplify the exposition.

We claim that $\Hol(G)$ is a symmetric smooth convolution algebra on~$G$.  It
is clear that the left regular representation of~$G$ on $\Hol(G)$ is smooth.
The difference between left and right convolution by~$D$ is given by the
adjoint representation of~$G$ on~$\LG$.  Since this has at most exponential
growth, the operator $\phi\mapsto\phi\ast D$ is a bounded linear operator on
$\Hol(G)$ for any $D\in\UE(\LG)$.  This means that the right regular
representation of~$G$ on $\Hol(G)$ is smooth as well.  Thus $\Hol(G)$ is a
smooth convolution algebra on~$G$.  Of course, we could equally well have
defined $\Hol(G)$ as the smoothening of the right regular representation.
Since it does not matter on which side we smoothen the representation, the
same argument as in the discrete case shows that $\Hol(G)$ is symmetric.

Now suppose that~$G$ is a connected nilpotent Lie group or, slightly more
generally, an extension of a compact group by a connected nilpotent Lie group.
Let $B(G)\defeq \bigcap_{k\in\N} \ell_1(G,(\ell+1)^k)$ and let $\Sch(G)$ be
the smoothening of the left regular representation of~$G$ on $B(G)$.  This
space can be described as above.  Since the adjoint representation of~$G$ has
polynomial growth in this case, the right and left regular representations
of~$G$ on $\Sch(G)$ have the same smoothening.  The same reasoning as above
shows that $\Sch(G)$ is a symmetric smooth convolution algebra on~$G$.
Finally, we define $B^\omega(G)\defeq \bigcup_{\alpha>1}
\ell_1(G,\alpha^\ell)$ and let $\Sch^\omega(G)$ be the smoothening of the left
regular representation on this space.  This is a symmetric smooth convolution
algebra for the same reasons.

We remark without proof that the smoothening of $\bigcap
\ell_p(G,\alpha^\ell)$ is equal to $\Hol(G)$ for all $p\in[1,\infty]$ and for
all almost connected locally compact groups.  Similarly, in the nilpotent case
the smoothenings of $\bigcap \ell_p(G,(\ell+1)^k)$ and $\bigcup
\ell_p(G,\alpha^\ell)$ are equal to $\Sch(G)$ and $\Sch^\omega(G)$,
respectively, for all $p\in[1,\infty]$.  This implies that the spaces
$\Hol(G)$, $\Sch(G)$, and $\Sch^\omega(G)$ are nuclear whenever they are
defined.  We do not prove this statement because we have no use for it here.

\begin{theorem}  \label{the:isocoh_nilpotent}
  The convolution algebras $\Sch(G)$, $\Sch^\omega(G)$ and $\Hol(G)$ are
  isocohomological if~$G$ is a connected nilpotent Lie group.
\end{theorem}

\begin{proof}
  We have already remarked above that these convolution algebras are
  symmetric.  The complex $\CCINF\Omega(G)$ is a projective resolution of the
  trivial representation of~$G$ by Proposition~\ref{pro:deRham_resolution}.
  As in the proof of Lemma~\ref{lem:basic_nilpotent}, one shows that
  $\T(G)\hot_G \CCINF\Omega(G)\cong\T\Omega(\LG)$.  This yields the assertion
  as above.
\end{proof}

Let~$G$ be a Lie group with finitely many connected components and let
$K\subseteq G$ be a maximal compact subgroup.  The general structure theory of
such groups asserts that $G/K$ is diffeomorphic to some~$\R^n$.  Hence
$\CCINF\Omega(G/K)\to\C$ is a projective $\CCINF(G)$\brd{}module resolution of
the trivial representation of~$G$ by Proposition~\ref{pro:deRham_resolution}.
We choose a $K$\nbd{}invariant inner product on the Lie algebra~$\LG$ of~$G$.
This defines a $G$\nbd{}invariant Riemannian metric on $G/K$.  Let~$\LGP$ be
the orthogonal complement of the Lie algebra of~$K$ inside~$\LG$.  We may
identify~$\LGP$ with the tangent space of $G/K$ at the identity coset~$K$.  As
above, unit speed geodesics in $G/K$ emanating from~$K$ are exactly the paths
of the form $\exp(tX)$ for some $X\in\LGP$ with $\norm{X}=1$.

\begin{theorem}  \label{the:isocoh_reductive}
  Suppose that the map $\exp\colon\LGP\to G/K$ described above is a
  diffeomorphism.  For instance, this is the case if~$\LG$ is a reductive Lie
  algebra and~$G$ and its centre have only finitely many connected components,
  or if~$G$ is an exponential Lie group.  Then $\Hol(G)$ is isocohomological.
\end{theorem}

\begin{proof}
  It is well-known that the homogeneous space $G/K$ is a CAT(0) space if~$\LG$
  is reductive and $G$ and its centre have only finitely many connected
  components.  This implies that the exponential map is a diffeomorphism
  (see~\cite{Bridson-Haefliger}).  Recall that~$G$ is exponential if
  $\exp\colon \LG\to G$ is a diffeomorphism.  This implies $K=\{1\}$ and hence
  that the map $\exp\colon \LGP\to G/K$ is a diffeomorphism.
  
  Since the exponential map $\LGP\to G/K$ is a diffeomorphism, it yields an
  isomorphism $\CCINF\Omega(G/K)\cong\CCINF\Omega(\LGP)$.  As above, the proof
  is finished by showing that this extends to an isomorphism between
  $\Hol(G)\hot_G \CCINF\Omega(G/K)$ and $\Hol\Omega(\LGP)$.  If $K=\{1\}$,
  that is, $G$ is an exponential Lie group, then the argument is almost
  literally the same as above.  The function $\sigma(\ad Z)$ above and its
  inverse automatically have exponential growth.  This is enough to prove that
  the different bases $(X_j)$ and $(Y_j)$ are equivalent for measuring the
  decay of derivatives, as long as we are only interested in spaces like
  $\Hol(G)$ which are closed under multiplication by functions of exponential
  growth.
  
  In general, we may proceed as follows.  The projection $G\to G/K$ induces a
  map $\Omega(G/K)\to\Omega(G)$.  We call a differential form
  $\omega\in\Omega(G)$ \emph{special} if it is in the range of this map.
  Equivalently, it is invariant under the right regular representation of~$K$
  and takes values in $\Lambda^*\LGP\subseteq \Lambda^*\LG$.  Let
  $X_1,\dots,X_m$ be a basis of~$\LG$ and view these elements as left
  invariant vector fields on~$G$.  We let $\Hol\Omega(G)$ be the space of all
  sections of $\Omega(G)$ whose coefficients with respect to the basis
  $dX_{i_1}\wedge\dotsb\wedge dX_{i_j}$ are in $\Hol(G)$.  Recall also that
  $f\in\Hol(G)$ if and only if $P(X_1,\dotsc,X_m)(f) = O(\alpha^\ell)$ for all
  $\alpha>1$ and all polynomials~$P$.  One can now check that $\Hol(G)\hot_G
  \CCINF\Omega(G/K)$ is isomorphic to the closed subspace of special forms in
  $\Hol\Omega(G)$.
  
  We may view $\LG\times K$ as a Lie group and then define
  $\Hol\Omega(\LG\times K)$ and the subspace of special differential forms as
  above.  Explicitly, special differential forms on $\LG\times K$ are defined
  using the coordinate projection $\LG\times K\to\LG$.  We use the length
  function $\ell(X,k)\defeq\norm{X}$ on $\LGP\times K$ and the basis
  $Y_1,\dotsc,Y_m$ of ${\LG\times K}$\brd{}invariant vector fields associated
  to the basis $(X_j)$ of~$\LG$.  One can show easily that $\Hol\Omega(\LGP)$
  is isomorphic to the subspace of special differential forms in
  $\Hol\Omega(\LG\times K)$.

  The hypotheses on~$G$ guarantee that the map
  $$
  \Theta\colon \LGP\times K\to G,
  \qquad (X,k)\mapsto \exp(X)\cdot k,
  $$
  is a diffeomorphism.  If $\omega\in\Omega(G)$, then~$\omega$ is special
  if and only if $\Theta^*(\omega)$ is special.  The function
  $\ell\circ\Theta^{-1}$ on~$G$ is equivalent to the length function that is
  used to define $\Hol(G)$.  The bases $(X_j)$ and $(Y_j)$ are related by
  multiplication by functions of at most exponential growth.  Therefore,
  $\Hol\Omega(G)\cong\Hol\Omega(\LG\times K)$, and the subspaces of special
  elements also coincide.  This means that $\Hol(G)\hot_G
  \CCINF\Omega(G/K)\cong\Hol\Omega(\LGP)$ as asserted.
\end{proof}

\section{Tempered crossed products and noncommutative tori}
\label{sec:crossed}

In this section, we consider crossed product algebras.  First we define them,
show that they are quasi-unital, and identify the essential modules with
covariant representations.  Then we find a sufficient criterion for
isocohomological embeddings in this context and apply it to noncommutative
tori.

\subsection{Crossed product algebras with compact support}
\label{sec:def_crossed}

Let~$G$ be a locally compact group and let~$B$ be a bornological algebra
equipped with a smooth representation of~$G$ by algebra automorphisms
$\beta\colon G\to\Aut(B)$.  The space $\CCINF(G,B)$ of smooth compactly
supported functions $G\to B$ can be defined most quickly as $\CCINF(G)\hot B$
(see also~\cite{Meyer:Smooth}).  Let~$d\lambda$ be a left invariant Haar
measure on~$G$.  We define the convolution on $\CCINF(G,B)$ by the usual
formula
\begin{equation}  \label{eq:crossed_convolution}
  \phi_1*\phi_2(g)
  \defeq
  \int_G \phi_1(h)\cdot \beta_h\phi_2(h^{-1}g) \,d\lambda(h).
\end{equation}
for all $\phi_1,\phi_2\in\CCINF(G,B)$.  This multiplication turns
$\CCINF(G,B)$ into a bornological algebra, which we denote $\CCINF(G)\cross
B$.

\begin{definition}  \label{def:covariant_rep}
  Suppose that~$B$ is essential.  An (essential) \emph{covariant
    representation of $G$ and~$B$} on a bornological vector space~$M$ is a
  pair $(\rho,\pi)$, where $\rho\colon B\to\End(M)$ is an essential
  $B$\nbd{}module structure on~$M$ and $\pi\colon G\to\Aut(M)$ is a smooth
  representation, such that $\pi(g)\rho(b)\pi(g^{-1})=\rho(\beta_gb)$ for all
  $g\in G$, $b\in B$.
\end{definition}

A covariant representation can be integrated to a module structure
$$
\CCINF(G)\cross B \to \End(M),
\qquad \phi\mapsto \int_G \rho(\phi(g)) \circ \pi(g)\,d\lambda(g).
$$
We want to show that all essential modules over $\CCINF(G)\cross B$ are of
this form.  As in the theory of $C^*$\nbd{}algebra crossed products, this is
done by embedding $B$ and~$G$ in the multiplier algebra $\Mult(\CCINF(G)\cross
B)$ by
\begin{alignat*}{2}
  (g\cdot\phi)(h) & \defeq \beta_g\phi(g^{-1}h),
  & \qquad (b\cdot \phi)(h) &\defeq b\phi(h),
  \\
  (\phi\cdot g)(h) &\defeq \phi(hg^{-1}) \mu_G(g)^{-1},
  & \qquad (\phi\cdot b)(h) &\defeq \phi(h) \beta_h(b),
\end{alignat*}
where~$\mu_G$ denotes the modular function of~$G$.  The representations of~$G$
on $\CCINF(G)\cross B$ by left and right multiplication are evidently smooth,
that is, $\CCINF(G)\cross B$ is an essential bimodule over $\CCINF(G)$.
If~$B$ is quasi-unital, then $\CCINF(G)\cross B$ is also an essential
$B$\nbd{}bimodule.  This is clear for the left module structure.  For the
right module structure, use the bounded linear operators $W,W^{-1}\colon
\CCINF(G,B)\to\CCINF(G,B)$,
\begin{equation}  \label{eq:multiplicative_invertible}
  W\phi(g)\defeq \beta_g\phi(g),
  \qquad
  W^{-1}\phi(g)\defeq \beta_g^{-1}\phi(g).
\end{equation}
They satisfy $W^{-1}(W(\phi)\cdot b)(h)=\phi(h)b$.  Hence $\CCINF(G)\cross
B\cong \CCINF(G)\hot B$ with the free module structure both as a left and a
right $B$\nbd{}module.

\begin{theorem}  \label{the:crossed_quasi-unital}
  If~$B$ is quasi-unital, then $\CCINF(G)\cross B$ is quasi-unital as well.
  The category of essential $\CCINF(G)\cross B$\brd{}modules is equivalent to
  the category of covariant representations of $G$ and~$B$.
\end{theorem}

\begin{proof}
  It is straightforward to construct approximate units in $\CCINF(G)\cross B$
  using approximate units in~$B$ and an approximation of~$\delta_1$ in
  $\CCINF(G)$ (see also~\cite{Meyer:Smooth} for a treatment of $\CCINF(G)$).
  Let $(\rho,\pi)$ be a covariant representation of $G$ and~$B$ on~$M$.  Since
  $\rho\colon B\to\End(M)$ is essential, we have a bounded linear section
  $\sigma_{BM}\colon M\to B\hot M$.  Define $\sigma\colon M\to\CCINF(G,B\hot
  M)\cong (\CCINF(G)\cross B)\hot M$ by
  $$
  \sigma(m)(g)\defeq \psi(g)\cdot (\ID_B\hot\pi)(g)^{-1} \sigma_{BM}(m)
  $$
  for some function $\psi\in\CCINF(G)$ with $\int_G
  \psi(g)\,d\lambda(g)=1$.  This is a bounded linear section for the action
  $(\CCINF(G)\cross B) \hot M\to M$.  Thus covariant representations yield
  essential modules.  Applying this to left and right module structure on
  $\CCINF(G)\cross B$, we see that $\CCINF(G)\cross B$ is quasi-unital.
  Conversely, if~$M$ is an essential $\CCINF(G)\cross B$\brd{}module, we get a
  smooth representation of~$G$ and an essential $B$\nbd{}module structure
  using the morphisms $\CCINF(G),B\dashrightarrow \CCINF(G)\cross B$.  They
  satisfy the covariance condition and their integrated form is the given
  $\CCINF(G)\cross B$\brd{}module structure.  One checks easily that the above
  constructions are functorial and inverse to each other.  Thus we have an
  isomorphism of categories.
\end{proof}

\subsection{Tempered crossed product algebras}
\label{sec:tempered_crossed}

Now let $\T(G)$ be a smooth convolution algebra on~$G$.  Let~$B$ be as above
and suppose, in addition, that the operators $W$ and~$W^{-1}$
in~\eqref{eq:multiplicative_invertible} extend to bounded linear operators on
$\T(G,B)\defeq \T(G)\hot B$.  Using this and that $\T(G)$ is a convolution
algebra, one can check easily that the
convolution~\eqref{eq:crossed_convolution} extends to a bounded bilinear map
on $\T(G,B)$.  We let $\T(G)\cross B$ be $\T(G,B)$ equipped with this
multiplication.

\begin{proposition}  \label{pro:tempered_cross}
  The algebra $\T(G)\cross B$ is quasi-unital and the embedding
  $$
  i\colon \CCINF(G)\cross B\to \T(G)\cross B
  $$
  is a proper morphism.  The category $\Mod(\T(G)\cross B)$ is isomorphic
  to the category of covariant representations $(\rho,\pi)$ for which
  $\pi\colon G\to\Aut(M)$ is $\T(G)$\brd{}tempered.
\end{proposition}

\begin{proof}
  The same argument as for $\CCINF(G)\cross B$ shows that $\T(G)\cross B$ is
  an essential $B$\nbd{}bimodule and that the representations of~$G$ on
  $\T(G)\cross B$ by left and right multiplication are smooth.  Hence
  $\T(G)\cross B$ carries a covariant representation of $G$ and~$B$.  By
  Theorem~\ref{the:crossed_quasi-unital}, $\T(G)\cross B$ is an essential
  bimodule over $\CCINF(G)\cross B$.  Lemma~\ref{lem:embed_quasiunital}
  implies that $\T(G)\cross B$ is quasi-unital and that the embedding of
  $\CCINF(G)\cross B$ is a proper morphism.  Since this homomorphism has dense
  range, the induced functor $\Mod((\T(G)\cross B)^+)\to \Mod((\CCINF(G)\cross
  B)^+)$ is fully faithful.  Proposition~\ref{pro:adjointness_proper_morphism}
  shows that the smoothening functor on $\Mod((\T(G)\cross B)^+)$ is the
  restriction of the smoothening functor on $\Mod((\CCINF(G)\cross B)^+)$.
  Thus a module over $\T(G)\cross B$ is essential if and only if it is
  essential as a module over $\CCINF(G)\cross B$.  By
  Theorem~\ref{the:crossed_quasi-unital}, $\Mod(\T(G)\cross B)$ is isomorphic
  to the category of covariant representations $(\pi,\rho)$ whose integrated
  form extends to a bounded homomorphism on $\T(G)\cross B$.  The latter just
  means that the integrated form of~$\pi$ extends to $\T(G)$.
\end{proof}

\begin{theorem}  \label{the:isocoh_crossed}
  Let~$G$ be a locally compact group and let $\T(G)$ be an isocohomological
  smooth convolution algebra on~$G$.  Let $B_1$ and~$B_2$ be quasi-unital
  bornological algebras, equipped with actions of~$G$ by automorphisms.
  Suppose that the operators $W$ and~$W^{-1}$
  in~\eqref{eq:multiplicative_invertible} extend to bounded operators on
  $\T(G)\hot B_2$, so that $\T(G)\cross B_2$ is defined.  Let $i\colon
  B_1\dashrightarrow B_2$ be a $G$\nbd{}equivariant isocohomological morphism.
  
  Then the induced morphism $\CCINF(G)\cross B_1\dashrightarrow \T(G)\cross
  B_2$ is isocohomological.
\end{theorem}

\begin{proof}
  By Proposition~\ref{pro:compose_isocoh}, it suffices to prove that the
  morphisms $\CCINF(G)\cross B_1\dashrightarrow\CCINF(G)\cross B_2$ and
  $\CCINF(G)\cross B_2\dashrightarrow\T(G)\cross B_2$ are isocohomological.
  
  The morphism $i^B\colon B\dashrightarrow \CCINF(G)\cross B$ gives rise to a
  functor
  $$
  i^B_!\colon \Mod(B)\to\Mod(\CCINF(G)\cross B),
  \qquad
  M\mapsto i^B_!(M) = (\CCINF(G)\cross B) \hot_B M.
  $$
  Since $\CCINF(G)\cross B$ is projective as a right $B$\nbd{}module, this
  functor is exact, that is, $i^B_!=\Left i^B_!$.  Moreover, $\CCINF(G)\cross
  B_2 \cong i^{B_1}_!(B_2) \cong (\CCINF(G)\cross B_1) \Lhot_{B_1} B_2$ and
  hence
  \begin{multline*}
    \CCINF(G)\cross B_2 \Lhot_{\CCINF(G)\cross B_1} \CCINF(G)\cross B_2
    \cong \CCINF(G)\cross B_2 \Lhot_{\CCINF(G)\cross B_1}
      \CCINF(G)\cross B_1 \Lhot_{B_1} B_2
    \\ \cong (\CCINF(G)\cross B_2) \Lhot_{B_1} B_2
    \cong \CCINF(G) \hot (B_2\Lhot_{B_1} B_2).
  \end{multline*}
  Thus $\CCINF(G)\cross B_1\dashrightarrow\CCINF(G)\cross B_2$ is
  isocohomological if and only if $B_1\dashrightarrow B_2$ is.
  
  Next we let $B=B_2$ and show that $\CCINF(G)\cross
  B\dashrightarrow\T(G)\cross B$ is isocohomological.
  The morphism $i^G\colon \CCINF(G)\dashrightarrow \CCINF(G)\cross B$ gives
  rise to a functor
  $$
  i^G_!\colon \Mod(G)\to\Mod(\CCINF(G)\cross B),
  \qquad
  M\mapsto i^G_!(M) = (\CCINF(G)\cross B) \hot_G M.
  $$
  Since $\CCINF(G)\cross B$ is projective as a right
  $\CCINF(G)$\nbd{}module, this functor is exact, that is, $i^G_!=\Left
  i^G_!$.  One checks easily that $\T(G)\cross B \cong i^G_! \T(G)$, using
  that $W$ and~$W^{-1}$ are bounded on $\T(G)\hot B$.  As above, this implies
  that $\CCINF(G)\cross B\to\T(G)\cross B$ is isocohomological if and only if
  $\T(G)$ is isocohomological.
\end{proof}

\subsection{Application to noncommutative tori}
\label{sec:app_tori}

In this section, we show how Alain Connes's computations for non-commutative
tori in~\cite{Connes:Noncommutative_Diffgeo} fit into our framework.  Our
method obviously extends to non-commutative tori of higher rank and to some
other crossed product algebras.  We do not pursue this here in order to
formulate results very concretely.  Let $\theta\in\R$ and let
$\Torus^2_\theta$ be the noncommutative $2$\nbd{}torus defined by invertible
generators $U,V$ satisfying the relation $UV=\exp(2\pi\ima\theta)VU$ for some
$\theta\in\R$.  We define bornological algebras
$$
\Pol(\Torus^2_\theta) \subseteq \Hol(\Torus^2_\theta)
\subseteq \Sch^\omega(\Torus^2_\theta) \subseteq \Sch(\Torus^2_\theta)
$$
of polynomial, holomorphic, real analytic, and smooth functions on
$\Torus^2_\theta$.  These spaces are equal to $\C[\Z^2]$, $\Hol(\Z^2)$,
$\Sch^\omega(\Z^2)$, and $\Sch(\Z^2)$ as bornological vector spaces, and
equipped with the product defined by $UV=\exp(2\pi\ima\theta)VU$.  For
$\theta=0$, we get the algebras of Laurent series, of holomorphic functions on
$(\C^\times)^2$, of real analytic functions on~$\Torus^2$, and of smooth
functions on $\Torus^2$, respectively.  Thus one may also use the notation
$\Pol((\C^\times)^2_\theta)$ and $\Hol((\C^\times)^2_\theta)$ to stress that
these algebras live on the complex torus $(\C^\times)^2$.

Let $\T(\Torus^2_\theta)$ be one of these bornological algebras.  We remark
that $\Pol(\Torus^2_\theta)$ carries the fine bornology,
$\Hol(\Torus^2_\theta)$ and $\Sch(\Torus^2_\theta)$ are Fréchet-Schwartz
spaces equipped with the precompact bornology, and
$\Sch^\omega(\Torus^2_\theta)$ is a Silva space
(see~\cite{Meyer:Combable_poly}).  However, the only property that we need is
that $\T(\Z)\hot\T(\Z)\cong\T(\Z^2)$.

Let~$\Z$ act on the convolution algebra $\T(\Z)$ by $n\cdot V^m\defeq
\exp(2\pi\ima\theta\cdot nm)V^m$.  This defines a representation of~$\Z$ by
automorphisms.  The operator~$W$ in~\eqref{eq:multiplicative_invertible} is
given by $W(U^mV^n) = \exp(2\pi\ima\theta\cdot mn) U^mV^n$.  This is evidently
a bounded linear operator on $\T(\Z\times\Z)=\T(\Z)\hot\T(\Z)$.  Hence we can
form the crossed product algebra $\T(\Z)\cross\T(\Z)$.  Of course, this is
nothing but $\T(\Torus^2_\theta)$.

\begin{theorem}  \label{the:isocoh_nctorus}
  The embeddings of $\Pol(\Torus^2_\theta)$ in $\Hol(\Torus^2_\theta)$,
  $\Sch^\omega(\Torus^2_\theta)$ and $\Sch(\Torus^2_\theta)$ are
  isocohomological.
\end{theorem}

\begin{proof}
  Theorem~\ref{the:isocoh_coarse} shows that the embedding $\C[\Z]\to\T(\Z)$
  is isocohomological.  Therefore, the embedding
  $\Pol(\Torus^2_\theta)=\C[\Z]\cross\C[\Z] \to
  \T(\Z)\cross\T(\Z)=\T(\Torus^2_\theta)$ is isocohomological by
  Theorem~\ref{the:isocoh_crossed}.
\end{proof}

Let~$K_\bullet$ be the complex of free $\Pol(\Torus^2_\theta)$\brd{}bimodules
$$
0 \to
\Pol(\Torus^2_\theta) \otimes \Pol(\Torus^2_\theta)
\overset{b_2}\longrightarrow
\bigl[\Pol(\Torus^2_\theta) \otimes \Pol(\Torus^2_\theta)\bigr]^2
\overset{b_1}\longrightarrow
\Pol(\Torus^2_\theta) \otimes \Pol(\Torus^2_\theta)
\to 0
$$
augmented by a map $K_0\overset{b_0}\to\Pol(\Torus^2_\theta)$, where
\begin{align*}
  b_0(x\otimes y) &\defeq x\cdot y,
  \\
  b_1(x_1\otimes y_1,x_2\otimes y_2)
  &\defeq x_1\otimes y_1 - x_1U^{-1}\otimes Uy_1
  + x_2\otimes y_2 - x_2V^{-1}\otimes Vy_2,
  \\
  b_2(x\otimes y) &\defeq
  (x\otimes y - xV^{-1}\otimes Vy, xU^{-1}\otimes Uy-x\otimes y).
\end{align*}
It is easy to check that $K_\bullet\to\Pol(\Torus^2_\theta)$ is a free
$\Pol(\Torus^2_\theta)$\brd{}bimodule resolution of $\Pol(\Torus^2_\theta)$.
This purely algebraic statement is not affected if we equip all spaces
in~$K_\bullet$ with the fine bornology.

Since $\Pol(\Torus^2_\theta)$ is free as a right module, the chain complex
above is still contractible as a complex of right modules.  Hence
$K_\bullet\hot_{\Pol(\Torus^2_\theta)} M\to M$ is a free resolution of~$M$ for
any $M\inOb\Mod(\Pol(\Torus^2_\theta))$.  As a result,
$\Right\Hom_{\Pol(\Torus^2_\theta)}(M,N)$ and
$M\Lhot_{\Pol(\Torus^2_\theta)} N$ can be identified with the
bornological chain complexes
\begin{gather*}
  0 \to \Hom_{\Pol(\Torus^2_\theta)}(M,N)
  \overset{b_2}\longrightarrow \Hom_{\Pol(\Torus^2_\theta)}(M,N)^2
  \overset{b_1}\longrightarrow \Hom_{\Pol(\Torus^2_\theta)}(M,N)
  \to 0
  \\
  0 \to M\hot N
  \overset{b_2}\longrightarrow (M\hot N)^2
  \overset{b_1}\longrightarrow M\hot N
  \to 0.
\end{gather*}
Theorem~\ref{the:isocoh_nctorus} implies that the same chain complexes also
compute $\Right\Hom_{\T(\Torus^2_\theta)}(M,N)$ and
$M\Lhot_{\T(\Torus^2_\theta)} N$.  In particular, $\T(\Torus^2_\theta)$ has
cohomological dimension~$2$.

The \emph{commutator quotient} $M/\natural$ of a bimodule~$M$ is the quotient
of~$M$ by the closed linear span of elements of the form $am-ma$ with $a\in
A$, $m\in M$.  The Hochschild homology of a quasi-unital algebra~$A$ is equal
to the homology of $P_\bullet/\natural$ for any projective bimodule resolution
$P_\bullet$ of~$A$ because such algebras are H-unital
(see~\cite{Loday:Cyclic}).  Theorem~\ref{the:isocoh}.1 provides such a
resolution for $\T(\Torus^2_\theta)$.  The resulting commutator quotient
complex is
$$
\T(\Torus^2_\theta)
\overset{\delta_2}\longrightarrow \T(\Torus^2_\theta)^2
\overset{\delta_1}\longrightarrow \T(\Torus^2_\theta)
$$
with $\delta_2(x)=(x-VxV^{-1},UxU^{-1}-x)$,
$\delta_1(x_1,x_2)=x_1-Ux_1U^{-1}+x_2-Vx_2V^{-1}$.  We can identify this chain
complex with $L\hot L$, where~$L$ is the chain complex
$$
\T(\Z) \overset{\alpha}\to\T(\Z),
\qquad \alpha f(m)\defeq (1-\exp(2\pi\ima\theta m))f(m),
$$
concentrated in degrees $0$ and~$1$.  The kernel of~$\alpha$ is spanned
by~$\delta_0$ if~$\theta$ is irrational.  The restriction of~$\alpha$ to
$\C[\Z^*]$ is invertible.  Thus the chain complex~$L$ for
$\Pol(\Torus^2_\theta)$ is homotopy equivalent to $\C\overset{0}\to\C$.  As a
result, the Hochschild homology of $\Pol(\Torus^2_\theta)$ is given by $\C$,
$\C^2$, and~$\C$ in dimensions $0$, $1$, and~$2$, and vanishes in higher
dimensions.

If we replace $\C[\Z^*]$ by $\T(\Z^*)$, the restriction of~$\alpha$ need not
be invertible any more.  The issue is whether or not the function $m\mapsto
(1-\exp(2\pi\ima\theta m))^{-1}$ on~$\Z^*$ has polynomial, subexponential, or
exponential growth, depending on whether we consider $\Sch(\Z)$,
$\Sch^\omega(\Z)$, or $\Hol(\Z)$.  It is observed by Alain Connes
(\cite{Connes:Noncommutative_Diffgeo}) that this question for $\Sch(\Z)$
depends on the Diophantine approximation properties of the parameter~$\theta$.
The same holds for the other function spaces.  The Hochschild homologies of
$\T(\Torus^2_\theta)$ and $\Pol(\Torus^2_\theta)$ agree if and only if
$(1-\exp(2\pi\ima\theta m))^{-1}$ satisfies the appropriate growth condition.
Otherwise, the Hochschild homology of $\T(\Torus^2_\theta)$ contains some
additional infinite dimensional non-Hausdorff spaces in dimensions $1$
and~$0$.  We remark that these additions disappear in periodic cyclic
homology.  The reason is the gauge action of~$\Torus^2$ on
$\T(\Torus^2_\theta)$.  By homotopy invariance, only the gauge invariant part
of the Hochschild homology contributes to the periodic cyclic homology.

\end{document}